\documentclass[12pt,a4paper]{article}
\usepackage{amssymb}
\usepackage{extarrows}
\usepackage{amsmath}
\usepackage{cite}
\newtheorem{theorem}{Theorem}[section]
\newtheorem{definition}[theorem]{Definition}
\newtheorem{lemma}[theorem]{Lemma}

\newtheorem{proposition}[theorem]{Proposition}
\newtheorem{example}[theorem]{Example}

\newtheorem{theorem(Composition-Diamond Lemma)}[theorem]{Theorem(Composition-Diamond Lemma)}

\begin{document}

\title{Gr\"{o}bner-Shirshov bases for associative conformal modules\footnote{Supported by the NNSF of China (no. 11571121) and the Science and Technology
Program of Guangzhou (no. 201707010137).}}
\author{
 Yuqun
Chen \\
{\small \ School of Mathematical Sciences, South China Normal
University}\\
{\small Guangzhou 510631, P. R. China}\\
{\small yqchen@scnu.edu.cn}\\
\\
Lili Ni\footnote{Supported by the talent fund of Taishan University (no. Y-01-2017001).}\\
{\small \ School of Mathematics and Statistics, Taishan
University}\\
{\small Taian 271000, P. R. China}\\
{\small nilili2009@163.com}}

\date{}
\maketitle

\maketitle \noindent\textbf{Abstract:} We construct  free modules
over  an
associative conformal algebra. We establish Composition-Diamond
lemma for associative conformal modules. As applications,
Gr\"{o}bner-Shirshov bases of the Virasoro conformal module and
module over the semidirect product of Virasoro conformal algebra and
current algebra are given respectively.


\noindent {\bf AMS} Mathematics Subject Classification(2000): 17B69,
16S15, 13P10

\section{Introduction}

The subject of conformal algebras is closely related to vertex
algebras (see,  V. Kac \cite{Kac96}). Implicitly, vertex algebras
were introduced by Belavin, Polyakov, and Zamolodchikov in 1984
\cite{Belavin}. Explicitly, the definition of vertex algebras was
given by R. Borcherds in 1986 \cite{Bor86}, which led to his
solution of the Conway-Norton conjecture in the theory of finite
simple groups \cite{Bor92, Frenkel88}.  As pointed out by Kac
\cite{Kac97, Kac96}, conformal and vertex algebras provide a
rigorous mathematical study of the ¡°locality axiom¡± which came
from Wightman¡¯s axioms of quantum field theory \cite{Wightman}. M.
Roitman studied free (Lie and associative) conformal and vertex
algebras in \cite{Ro99}. Free vertex algebras were mentioned in the
original paper of Borcherds \cite{Bor86}. Since conformal and vertex
algebras are not varieties in the sense of universal algebra (see,
P.M. Cohn \cite{Cohn81}), the existence of free conformal and free
vertex algebras is not guaranteed by the general theory and should
be proved. It was done by Roitman  \cite{Ro99}. The free associative
conformal algebra generated by a set $B$ with a  locality function
$N(-,-):\ B\times B\rightarrow \mathbb{Z}_{\geq 0}$ is constructed
by Bokut, Fong, and Ke in 1997 \cite{BFK0}.

Conformal
module is a basic tool for the construction of free field realization of infinite dimensional
Lie (super)algebras in conformal field theory.  Finite irreducible conformal modules over the Virasoro
conformal algebra were determined in \cite{SJVK}.  The Lie conformal algebra of a Block type was introduced and
free intermediate series modules were classified in \cite{GaoXuYue}.

  Gr\"{o}bner bases and Gr\"{o}bner-Shirshov bases  were invented
independently by A.I. Shirshov for ideals of free (commutative,
anti-commutative) non-associative algebras {\cite{Shir3,Sh62a}}, free
Lie algebras \cite{Shir3} and implicitly free associative
algebras \cite{Shir3}  (see also \cite{Be78,Bo76}), by H.
Hironaka \cite{Hi64} for ideals of the power series algebras (both
formal and convergent), and by B. Buchberger \cite{Bu70} for ideals
of the polynomial algebras.

Gr\"{o}bner bases and Gr\"{o}bner-Shirshov bases theories have been
proved to be very useful in different branches of mathematics,
including commutative algebra and combinatorial algebra. It is a
powerful tool to solve the following classical problems: normal
form; word problem; conjugacy problem; rewriting system; automaton;
embedding theorem; PBW theorem;  extension; homology;  growth
function; Dehn function; complexity; etc.

Up to now, different versions of Composition-Diamond lemma are known
for the following classes of algebras apart those mentioned above:
(color) Lie super-algebras \cite{Mik89, Mik92,Mik96}, tensor product
of a free algebra and a polynomial algebra \cite{MZ}, tensor product
of two free algebras \cite{BCC08}, Lie $p$-algebras \cite{Mik92},
associative conformal algebras \cite{BFK04, nll},  shuffle operads
\cite{DK10}, modules \cite{Chi,KL} (see also \cite{CCZ}),
right-symmetric algebras \cite{BCLi08}, dialgebras \cite{BCL08},
associative algebras with multiple operators \cite{BCQ08},
Rota-Baxter algebras \cite{BCD08}, Lie algebras over a polynomial
algebra \cite{BCC11},  matabelian Lie algebras \cite{CC12},
semirings \cite{BCM13}, and so on.

A Composition-Diamond lemma for
associative conformal algebras is firstly established
by  Bokut, Fong, and Ke in 2004  \cite{BFK04} which
claims that if (i) $S$ is a Gr\"{o}bner-Shirshov basis
in $C(B,N)$, then (ii) the set of $S$-irreducible words
is a linear basis of the quotient conformal algebra
$C(B,N|S)$, but not conversely.
By introducing some new definitions of normal $S$-words,
compositions and compositions to be trivial, a
new Composition-Diamond lemma for associative conformal
algebras is established by  Ni and Chen in 2016  \cite{nll}
which makes the conditions  (i) and (ii) equivalent.

In this article we establish  Gr\"{o}bner-Shirshov bases method for
modules over an associative conformal algebra and give some
applications.

 The article is organized as follows. In section 2, we introduce the concepts of conformal
 algebra, associative (Lie) conformal algebra, and modules over an associative (Lie) conformal algebra.
 In section 3, we construct  free modules over a free
 associative conformal algebra and free modules over an  associative conformal algebra.
 In section 4, we establish  Composition-Diamond lemma for associative conformal modules.
 In section 5, we give some applications of the Composition-Diamond lemma for associative
 conformal modules: Gr\"{o}bner-Shirshov bases of the Virasoro conformal module and module
 over the semidirect product of Virasoro conformal algebra and current algebra are given and
 then linear bases of them are obtained respectively.

\section{Preliminaries}

In this section, we introduce some related concepts.

\subsection{Associtive (Lie) conformal algebras}

We begin with the formal definition of a conformal algebra.

\begin{definition}(\cite{BFK0,BoFKK00,Ro99,Kac96,BFK00})
A conformal algebra $C=(C,(n),n\in\mathbb{Z}_{\geq 0},D)$ is a
linear space over a field $\mathbf{k}$ of characteristic $0$,
equipped with bilinear multiplications $a_{(n)}b$,
$n\in\mathbb{Z}_{\geq 0}=\{0,1,2,\cdots \}$,
and a linear map $D$, such that the following axioms are valid:
\begin{enumerate}
\item[(C1)]\ (locality) For any $a, b\in C$, there exists a nonnegative integer $N(a,b)$ such that $a_{(n)}b=0$ for $n\geq N(a,b)$\ ($N(a, b)$ is called the order of locality of $a$ and $b$);
\item[(C2)]\ $D(a_{(n)}b)=Da_{(n)}b+a_{(n)}Db$ for any $a,b\in C$
and $n\in\mathbb{Z}_{\geq 0}$;
\item[(C3)]\ $Da_{(n)}b=-na_{(n-1)}b$ for any $a,b\in C$ and
$n\in\mathbb{Z}_{\geq 0}$, and $Da_{(0)}b=0$.
\end{enumerate}

A conformal algebra $C$ is called  {\it associative} if the following identity holds for
all $a,b,c\in C,\ m,n\in\mathbb{Z}_{\geq 0}$,
$$
(a_{(n)}b)_{(m)}c=\sum_{t\geq 0}(-1)^{t}\binom{n}{t}a_{(n-t)}(b_{(m+t)}c).
$$

A conformal algebra $L=\langle L, [n], n\in\mathbb{Z}_{\geq 0},
D\rangle$  is called a Lie conformal algebra if $L$ satisfies the
following two axioms:

$\bullet$ (Anti-commutativity)\ $a_{[n]}b=-\{b_{[n]}a\}$, where
$$
\{b_{[n]}a\}=\sum\limits_{k\geq
0}(-1)^{n+k}\frac{1}{k!}D^{k}(b_{[n+k]}a).
$$

$\bullet$ (Jacobi identity)
$$
(a_{[n]}b)_{[m]}c=\sum\limits_{k\geq 0}(-1)^{k}
\binom{n}{k}(a_{[n-k]}(b_{[m+k]}c)-b_{[m+k]}(a_{[n-k]}c)).
$$
\end{definition}

\subsection{Modules over a Lie or associative conformal algebra}

\begin{definition}(\cite{SJVK,Kac96,ARe})
Let $C$ be an associative (Lie, resp.) conformal algebra.  An associative (Lie, resp.) conformal module $_CM$ is a $\mathbf{k}[D]$-module $M$
 endowed with a series of operations $(n):$ $C\times M\rightarrow M$,
 $n\in\mathbb{Z}_{\geq 0}$, such that for any $a,b\in C$ and $v\in M$,
 \begin{enumerate}
 \item[(i)]\ (locality) there exists a nonnegative integer $N(a,v)$ such that $a_{(n)}v=0$ for $n\geq N(a,v)$\ ($N(a, v)$ is called the order of locality of $a$ and $v$);
 \item[(ii)]\ $D(a_{(n)}v)=Da_{(n)}v+a_{(n)}Dv$ for $n\in\mathbb{Z}_{\geq 0}$;
 \item[(iii)]\ $Da_{(n)}v=-na_{(n-1)}v$ for $n\in\mathbb{Z}_{\geq 0}$, and $Da_{(0)}v=0$;
 \item[(iv)]\ $(a_{(n)}b)_{(m)}v=\sum\limits_{t\geq 0}(-1)^{t}\binom{n}{t}a_{(n-t)}(b_{(m+t)}v)$ if $C$ is associative;\\
$(a_{[n]}b)_{(m)}v=\sum\limits_{t\geq 0}(-1)^{t}\binom{n}{t}(a_{(n-t)}(b_{(m+t)}v)-b_{(m+t)}(a_{(n-t)}v))$ if $C$ is Lie.
 \end{enumerate}
\end{definition}

A  {\it submodule} $M_1$ of a conformal module
$_CM$ is a $\mathbf{k}[D]$-submodule
such that for all $n\in\mathbb{Z}_{\geq 0}$, $C_{(n)}M_1\subseteq M_1$. If $S\subseteq  _C\!M$, $subm(S)$ means
the submodule of $_CM$ generated by $S$.

\section{Free associative conformal modules}

\begin{definition}
An associative conformal module $mod_{_C}(Y)$ over an associative conformal algebra $C$
is called the free associative conformal module
generated by a set $Y$ with the locality function $N:C\times
Y\longrightarrow \mathbb{Z}_{\geq 0}$,  if for any associative conformal module $_CM$ and any mapping $\varepsilon :
Y\longrightarrow M$ with
$c_{(n)}\varepsilon (y)=0$ for all $c\in C,\ y\in Y$ and
$n\geq N(c,y)$, there exists a unique $C$-module homomorphism
$\varphi: mod_{_C}(Y)\rightarrow M$ such that the following diagram is commutative:

\setlength{\unitlength}{1cm}
\begin{picture}(7, 3)
\put(4.2,2.3){\vector(1,0){1.7}} \put(4.1, 2.0){\vector(0,-1){1.3}}
\put(5.9,2.1){\vector(-1,-1){1.6}} \put(3.9,0.2){$M$}
\put(6,2.2){$mod_{_C}(Y)$} \put(3.9,2.2){$Y$} \put(4.9, 2.4){$i$}
\put(3.7,1.3){$\varepsilon$} \put(5.3,1){$\exists!\varphi$}
\end{picture}

\noindent where $i$ is the inclusion map.
\end{definition}

In this section, we construct the free associative conformal module
generated by a set $Y$.

\subsection{Double free associative conformal modules}

Let $C(B,N)$  be the free associative conformal algebra generated by
$B$ with the locality function $N(-,-):\ B\times B\longrightarrow\mathbb{Z}_{\geq 0}$ over a field $\mathbf{k}$ of characteristic $0$. Let $Y$
be a set. We construct the free module over $C(B,N)$ generated by
$Y$ which is called double free associative conformal module.

We extend the mapping $N(-,-)$ to $B\times (B\cup Y)\rightarrow \mathbb{Z}_{\geq 0}$. Let
\begin{eqnarray*}
U&=&\{b_{1(n_{1})}b_{2(n_{2})}\cdots b_{k(n_{k})}D^{i}y\mid  b_{j},b_k\in B,\ 0\leq n_{j}<N(b_j, b_{j+1}), 1\leq j\leq k-1, \\
&& \ \ \ \ \ \ \ \ \ \ \ \ \ \ \ \ \ \ \ \ \ \ \ \ \ \ \ \ \ \ \ \ \ \ \ \ 0\leq n_k< N(b_k, y),\ i\geq 0,\ k\geq
0,\ y\in Y\},\\
T&=&\{b_{1(n_{1})}\cdots b_{k(n_{k})}D^{i}b_{k+1}\mid
b_{j},b_{k+1}\in B, 0\leq n_{j}<N(b_j, b_{j+1}), 1\leq j\leq k, i, k\geq 0\}
\end{eqnarray*}
and
$$
[u]:=b_{1(n_{1})}(b_{2(n_{2})}(\cdots (b_{k(n_{k})}D^{i}y)\cdots)).
$$
If $k=0$, then $[u]=D^{i}y$.
We call $k+1$ the length of $u$, denoted by  $|u|$.
Note that $[\ ]$ means right normed bracketing.

Denote  $span_{\mathbf{k}}[U]$ the $\mathbf{k}$-linear space with a $\mathbf{k}$-basis $[U]:=\{[u]\mid u\in U\}$. We will make $span_{\mathbf{k}}[U]$ to be a $C(B,N)$-module.

\ \

Note that the set $[T]:=\{[a]\mid a\in T\}$ a $\mathbf{k}$-basis of  the free associative conformal algebra $C(B,N)$.

We define a scheme of {\it algorithm} for $C(B,N)$-module
$span_{\mathbf{k}}[U]$.
For any $a\in T,\ b\in B,\ y\in Y,\ n\geq 0,\ u\in U$, we define $[a]_{(n)}[u]$ as follows.

\begin{enumerate}
 \item[(i)]\ $
b_{(n)}[u]=\left\{
\begin{aligned}
& 0, \ \ \ \ \ \ \ \ \ \ \ \ \ \ \  \ \ \ \ \ \ \ \ \ \  \ \ \ \ \ \ \ \ \ \ \ \ \ \ \ \ \ \ \ \ \ \mbox{if} \ [u]=y,\ n\geq N(b,y) ,    \\
&-\sum\limits_{t\geq 1}(-1)^{t}\binom{i}{t}\frac{n!}{(n-t)!}b_{(n-t)}D^{i-t}y,
\ \ \mbox{if} \ [u]=D^{i}y, \ i> 0,\ n\geq N(b, y),  \\
&-\sum\limits_{t\geq 1}(-1)^{t}\binom{n}{t}b_{(n-t)}(b_{1(m+t)}[u_1]) ,
\ \ \ \ \mbox{if} \ [u]=b_{1(m)}[u_{1}],\ n\geq N(b, b_1).  \\
\end{aligned}
\right.
$

\item[(ii)]\ $
[a]_{(n)}[u]=\left\{
\begin{aligned}
& 0,  \ \ \ \ \ \ \ \ \ \ \ \ \ \ \ \ \ \ \ \ \ \ \ \ \ \ \ \ \ \ \ \ \ \ \ \ \ \ \ \ \ \mbox{if} \ \ [a]=D^{i}b, \ i>n,   \\
&(-1)^{i}\frac{n!}{(n-i)!}b_{(n-i)}[u], \ \ \ \ \ \ \ \ \ \ \ \ \ \ \   \mbox{if} \ \ [a]=D^{i}b, \ 0<i\leq n,  \\
&\sum_{t\geq 0}(-1)^{t}\binom{m}{t}b_{(m-t)}([c]_{(n+t)}[u]),
\ \  \mbox{if} \ \ [a]=b_{(m)}[c]. \\
\end{aligned}
\right.
$

\item[(iii)]\ For $u=b_{1(n_{1})}b_{2(n_{2})}\cdots b_{k(n_{k})}D^{i}y\in U$,
\begin{align*}
D([u])=
&\sum\limits_{j=1}^{k}[b_{1(n_{1})}\cdots
b_{j-1(n_{j-1})}Db_{j(n_{j})}b_{j+1(n_{j+1})}\cdots
 b_{k(n_{k})}D^{i}y] \ \ \ \ \ \ \ \ \  \\
&+[b_{1(n_{1})}b_{2(n_{2})}\cdots b_{k(n_{k})}D^{i+1}y].
\end{align*}
\end{enumerate}

\begin{lemma}\label{l2.0}
For any $b\in B,\ y\in Y$, $n\geq 0, \ i\geq 0$, we have
\begin{equation}\label{e2.2}
D(b_{(n)}D^{i}y)=Db_{(n)}D^{i}y+b_{(n)}D^{i+1}y
\end{equation}
\end{lemma}
\noindent{\bf Proof.} If $0\leq n<N(b, y)$ or $i=1$, the result follows from the definition.
We may assume that $n\geq N(b, y)$ and $i>1$.

Induction on $n$. Let $N:=N(b, y)$. The left hand side of (\ref{e2.2}) is
\begin{align*}
D(b_{(N)}D^{i}y)
&=D(-\sum\limits_{t\geq 1}(-1)^{t}\binom{i}{t}\frac{N!}{(N-t)!}b_{(N-t)}D^{i-t}y) \\
&=-\sum\limits_{t\geq 1}(-1)^{t}\binom{i}{t}\frac{N!}{(N-t)!}D(b_{(N-t)}D^{i-t}y) \\
&=A_1+A_2,
\end{align*}
where
$
A_1:=\sum\limits_{t\geq 1}(-1)^{t}\binom{i}{t}\frac{N!}{(N-t-1)!}b_{(N-t-1)}D^{i-t}y
=\sum\limits_{p\geq 2}(-1)^{p-1}\binom{i}{p-1}\frac{N!}{(N-p)!}b_{(N-p)}D^{i-p+1}y
$
and $A_2:=-\sum\limits_{t\geq 1}(-1)^{t}\binom{i}{t}\frac{N!}{(N-t)!}b_{(N-t)}D^{i-t+1}y$.
Hence
\begin{align*}
A_1+A_2
&=-\sum\limits_{t\geq 2}(-1)^{t}\binom{i+1}{t}\frac{N!}{(N-t)!}b_{(N-t)}D^{i-t+1}y
+iNb_{(N-1)}D^{i}y \\
&=-\sum\limits_{t\geq 1}(-1)^{t}\binom{i+1}{t}\frac{N!}{(N-t)!}b_{(N-t)}D^{i-t+1}y
-Nb_{(N-1)}D^{i}y.
\end{align*}

The right hand side of (\ref{e2.2}) is
$$
Db_{(N)}D^{i}y+b_{(N)}D^{i+1}y
=-Nb_{(N-1)}D^{i}y-\sum\limits_{t\geq 1}(-1)^{t}
\binom{i+1}{t}\frac{N!}{(N-t)!}b_{(N-t)}D^{i-t+1}y .
$$
So the result holds for $n=N$. Assume that $n>N$. Then the left hand side
of  (\ref{e2.2}) is
$$
D(b_{(n)}D^{i}y)
=-\sum\limits_{t\geq 1}(-1)^{t}\binom{i}{t}\frac{n!}{(n-t)!}D(b_{(n-t)}D^{i-t}y)
$$
and the right hand side of  (\ref{e2.2}) is
$$
Db_{(n)}D^{i}y+b_{(n)}D^{i+1}y
=-nb_{(n-1)}D^{i}y-\sum\limits_{t\geq 1}(-1)^{t}
\binom{i+1}{t}\frac{n!}{(n-t)!}b_{(n-t)}D^{i-t+1}y .
$$
By induction, we complete the proof.  \ \ \ \  $\square$

\begin{lemma}\label{l2.1}
For any $a\in T,\ u\in U$, there exists $N(a,[u])\in\mathbb{Z}_{\geq 0}$ such that $[a]_{(n)}[u]=0$ for all $n\geq N(a,[u])$.
\end{lemma}
\noindent{\bf Proof.}
Induction on $|a|$. If $|a|=1$, say, $[a]=b$, then $[a]_{(n)}[u]=b_{(n)}[u]$.
Now, we find, by induction on $|u|$, $N(b,[u])\in\mathbb{Z}_{\geq 0}$ such that $b_{(n)}[u]=0$ for all $n\geq N(b,[u])$.
Let $[u]=D^{i}y$.
Then by Lemma \ref{l2.0},
$$
b_{(n)}[u]=b_{(n)}D^{i}y=\sum\limits_{t\geq 0}\binom{n}{t}\frac{i!}{(i-t)!}D^{i-t}(b_{(n-t)}y).
$$
So $b_{(n)}D^{i}y=0$ if $n\geq N(b,y)+i$.
Hence, take $N(b,D^{i}y)=N(b,y)+i$ as required.

Assume that $|u|>1$ and $[u]=b_{1(m)}[u_1]$. Let
$N_0=N(b,b_1)+N(b_1, [u_1])-m-1$. Therefore, for $n\geq N_0\geq N(b,b_1)$, we have
$$
b_{(n)}(b_{1(m)}[u_1])=-\sum\limits_{t_1\geq 1}(-1)^{t_1}
\binom{n}{t_1}b_{(n-t_1)}(b_{1(m+t_1)}[u_1]).
$$
Since
\begin{align*}
m+t_1<N(b_1, [u_1])
&\Leftrightarrow t_1<N(b_1, [u_1])-m\Leftrightarrow n-t_1>n-N(b_1, [u_1])+m  \\
&\Leftrightarrow n-t_1>n-N_0+N(b,b_1)-1\Leftrightarrow n-t_1\geq N(b,b_1),
\end{align*}
we have
\begin{align*}
b_{(n)}(b_{1(m)}[u_1])
&=-\sum\limits_{t_1\geq 1}(-1)^{t_1}\binom{n}{t_1}b_{(n-t_1)}(b_{1(m+t_1)}[u_1]) \\
&=\sum\limits_{t_1,t_2\geq 1}(-1)^{t_1+t_2}\binom{n}{t_1}\binom{n-t_1}{t_2}
b_{(n-t_1-t_2)}(b_{1(m+t_1+t_2)}[u_1]) \\
&=\cdots   \\
&=\sum\limits_{t_j\geq 1}(-1)^{t+k}\binom{n}{t_1}\binom{n-t_1}{t_2}
\cdots \binom{n-t+t_{k}}{t_k}b_{(n-t)}(b_{1(m+t)}[u_1]),
\end{align*}
where $t=t_1+\cdots+t_k,\ n-t<N(b,b_1)$. Thus $m+t\geq N(b_1,[u_1])$.
So $b_{(n)}[u]=0$ for each $n\geq N_0$.

Assume that $|a|>1$ and $[a]=b_{(m)}[c]$. Since
$$
(b_{(m)}[c])_{(n)}[u]=\sum_{p\geq 0}(-1)^{p}\binom{m}{p}b_{(m-p)}([c]_{(n+p)}[u]),
$$
we can assume $n\geq N([c],[u])$. Therefore, $(b_{(m)}[c])_{(n)}[u]=0$. \ \ \ \  $\square$

\begin{lemma}\label{l2.2}
For any $a\in T,\ u\in U,\ m\geq 0$, we have
\begin{equation}\label{e2.1}
D[a]_{(m)}[u]=-m[a]_{(m-1)}[u]
\end{equation}
\end{lemma}
\noindent{\bf Proof.}
Assume first $|a|=1$ and $[a]=D^{i}b$. By definition, the left hand side
of (\ref{e2.1}) is
$$
D[a]_{(m)}[u]=D^{i+1}b_{(m)}[u]=\left\{
\begin{aligned}
& 0,\ \ \ \ \ \ \ \ \ \ \ \ \ \ \ \ \ \ \ \ \ \ \ \ \ \ \ \ \ \ \ \ \ \ \ \ \ \ \ i+1>m,   \\
&(-1)^{i+1}\frac{m!}{(m-i-1)!}b_{(m-i-1)}[u],\ i+1\leq m.
\end{aligned}
\right.
$$
While the right hand side of (\ref{e2.1}) is
$$
-m[a]_{(m-1)}[u]=(-m)D^{i}b_{(m-1)}[u]=\left\{
\begin{aligned}
& 0,\ \ \ \ \ \ \ \ \ \ \ \ \ \ \ \ \ \ \ \ \ \ \ \ \ \ \ \ \ \ \ \ \ \ \ \ \ \ \ \ \ \  i>m-1,   \\
&(-1)^{i+1}m\frac{(m-1)!}{(m-i-1)!}b_{(m-i-1)}[u],\ i\leq m-1.
\end{aligned}
\right.
$$
So $D(D^{i}b)_{(m)}[u]=-m(D^{i}b)_{(m-1)}[u]$.

Assume $|a|>1$, $[a]=b_{(n)}[c]$. Then
\begin{eqnarray*}
&&D[a]_{(m)}[u]
=D(b_{(n)}[c])_{(m)}[u]=-n(b_{(n-1)}[c])_{(m)}[u]+(b_{(n)}D[c])_{(m)}[u]  \\
&=&-\sum_{t\geq 0}(-1)^{t}n\binom{n-1}{t}b_{(n-1-t)}([c]_{(m+t)}[u])
+\sum_{t\geq 0}(-1)^{t}\binom{n}{t}b_{(n-t)}(D[c]_{(m+t)}[u])\\
&=&\sum_{p\geq 1}(-1)^{p}n\binom{n-1}{p-1}b_{(n-p)}([c]_{(m+p-1)}[u])
-\sum_{t\geq 0}(-1)^{t}(m+t)\binom{n}{t}b_{(n-t)}([c]_{(m+t-1)}[u]) \\
&=&\sum_{t\geq 1}(-1)^{t}t\binom{n}{t}b_{(n-t)}([c]_{(m+t-1)}[u])
-\sum_{t\geq 0}(-1)^{t}(m+t)\binom{n}{t}b_{(n-t)}([c]_{(m+t-1)}[u]) \\
&=&-m\sum_{t\geq 0}(-1)^{t}\binom{n}{t}b_{(n-t)}([c]_{(m+t-1)}[u]) \\
&=&-m(b_{(n)}[c])_{(m-1)}[u]=-m[a]_{(m-1)}[u].
\end{eqnarray*}
Hence (\ref{e2.1}) is true. \ \ \ \  $\square$

\begin{lemma}\label{l2.3}
For any $a\in T,\ u\in U,\ m\geq 0$, we have
\begin{equation}\label{e2.0}
D([a]_{(m)}[u])=D[a]_{(m)}[u]+[a]_{(m)}D[u]
\end{equation}
\end{lemma}
\noindent{\bf Proof.}
Case 1. Let $[a]=b,\ [u]=D^{i}y$. By Lemma \ref{l2.0}, the identity (\ref{e2.0}) holds
for $m\geq 0$, $i\geq 0$.

Case 2. Let $[a]=b$ and $[u]=b_{1(n)}D^{i}y$. If $0\leq m< N(b,b_1)$,
the identity (\ref{e2.0}) is true by the definition.
Assume that $[u]=b_{1(n)}D^{i}y,\ m\geq N(b,b_1)$. Induction on $m$. Let $N:=N(b,b_1)$.
Then
$$
b_{(N)}(b_{1(n)}D^{i}y)=-\sum_{t\geq 1}(-1)^{t}\binom{N}{t}b_{(N-t)}(b_{1(n+t)}D^{i}y).
$$
So the left hand side of (\ref{e2.0}) is
\begin{align*}
D(b_{(N)}[u])
&=-\sum_{t\geq 1}(-1)^{t}\binom{N}{t}Db_{(N-t)}(b_{1(n+t)}D^{i}y)
-\sum_{t\geq 1}(-1)^{t}\binom{N}{t}b_{(N-t)}D(b_{1(n+t)}D^{i}y) \\
&=A_1+A_2+A_3,
\end{align*}
where
\begin{align*}
A_1&:=\sum_{t\geq 1}(-1)^{t}(N-t)\binom{N}{t}b_{(N-t-1)}(b_{1(n+t)}D^{i}y), \\
A_2&:=\sum_{t\geq 1}(-1)^{t}(n+t)\binom{N}{t}b_{(N-t)}(b_{1(n+t-1)}D^{i}y), \\
A_3&:=-\sum_{t\geq 1}(-1)^{t}\binom{N}{t}b_{(N-t)}(b_{1(n+t)}D^{i+1}y)
=b_{(N)}(b_{1(n)}D^{i+1}y).
\end{align*}
The right hand side of (\ref{e2.0}) is
\begin{align*}
Db_{(N)}[u]+b_{(N)}D[u]
&=-Nb_{(N-1)}(b_{1(n)}D^{i}y)+b_{(N)}(Db_{1(n)}D^{i}y)+b_{(N)}(b_{1(n)}D^{i+1}y) \\
&=-Nb_{(N-1)}(b_{1(n)}D^{i}y)+\sum_{t\geq 1}(-1)^{t}n\binom{N}{t}b_{(N-t)}(b_{1(n+t-1)}D^{i}y)+A_3.
\end{align*}
Since
\begin{align*}
A_1+A_2
&=\sum_{t\geq 1}(-1)^{t}(N-t)\binom{N}{t}b_{(N-t-1)}(b_{1(n+t)}D^{i}y) \\
&\ \ \ +\sum_{t\geq 1}(-1)^{t}(n+t)\binom{N}{t}b_{(N-t)}(b_{1(n+t-1)}D^{i}y) \\
&=\sum_{p\geq 2}(-1)^{p-1}(N-p+1)\binom{N}{p-1}b_{(N-p)}(b_{1(n+p-1)}D^{i}y) \\
&\ \ \ +\sum_{t\geq 1}(-1)^{t}(n+t)\binom{N}{t}b_{(N-t)}(b_{1(n+s-1)}D^{i}y) \\
&=\sum_{t\geq 1}(-1)^{t-1}t\binom{N}{t}b_{(N-t)}(b_{1(n+t-1)}D^{i}y)
-mNb_{(N-1)}(b_{1(n)}D^{i}y) \\
&\ \ \ +\sum_{t\geq 1}(-1)^{t}n\binom{N}{t}b_{(N-t)}(b_{1(n+t-1)}D^{i}y)
+\sum_{t\geq 1}(-1)^{t}t\binom{N}{t}b_{(N-t)}(b_{1(n+t-1)}D^{i}y)  \\
&=-Nb_{(N-1)}(b_{1(n)}D^{i}y)+
\sum_{t\geq 1}(-1)^{t}n\binom{N}{t}b_{(N-t)}(b_{1(n+t-1)}D^{i}y),
\end{align*}
we have $D(b_{(N)}[u])=Db_{(N)}[u]+b_{(N)}D[u]$. By induction on $m$, we can
get $D(b_{(m)}[u])=Db_{(m)}[u]+b_{(m)}D[u]$ for any $|u|=2,\ m\geq 0$.

Next, we use induction on $|u|$. Assume that $|u|>2$ and the result is true for $|u|<l$.
Let $|u|=l$. Then, we can repeat the argument by induction on $m$, and get the identity
$D(b_{(m)}[u])=Db_{(m)}[u]+b_{(m)}D[u]$ for any $[u]\in U,\ m\geq 0$. Hence,
the identity (\ref{e2.0}) is true when $|a|=1, |u|\geq 1$.

Case 3. Suppose $|a|\geq 1$, $[a]=b_{(p)}[c]$. By induction on $|a|+|u|$, the
left hand side of (\ref{e2.0}) is
\begin{align*}
D([a]_{(m)}[u])&=D((b_{(p)}[c])_{(m)}[u]) \\
&=\sum_{t\geq 0}(-1)^{t}\binom{p}{t}D(b_{(p-t)}([c]_{(m+t)}[u]))
=A_1+A_2+A_3,
\end{align*}
where
\begin{align*}
A_1&:=-\sum_{t\geq 0}(-1)^{t}(p-t)\binom{p}{t}b_{(p-t-1)}([c]_{(m+t)}[u]), \\
A_2&:=-\sum_{t\geq 0}(-1)^{t}(m+t)\binom{p}{t}b_{(p-t)}([c]_{(m+t-1)}[u]), \\
A_3&:=\sum_{t\geq 0}(-1)^{t}\binom{p}{t}b_{(p-t)}([c]_{(m+t)}D[u])
=(b_{(p)}[c])_{(m)}D[u]=[a]_{(m)}D[u].
\end{align*}
Since
\begin{eqnarray*}
&&D[a]_{(m)}[u]=D(b_{(p)}[c])_{(m)}[u]=(Db_{(p)}[c])_{(m)}[u]+(b_{(p)}D[c])_{(m)}[u]\\
&=&-\sum_{t\geq 0}(-1)^{t}p\binom{p-1}{t}b_{(p-t-1)}([c]_{(m+t)}[u])
-\sum_{t\geq 0}(-1)^{t}(m+t)\binom{p}{t}b_{(p-t)}([c]_{(m-1+t)}[u]) \\
&=&A_1+A_2,
\end{eqnarray*}
the identity (\ref{e2.0}) is true. \ \ \ \  $\square$

\begin{lemma}\label{l2.5}
Let $n,m\geq 0,\ a, c\in T,\ u\in U$. Then
\begin{eqnarray}\label{e2.3}
([a]_{(n)}[c])_{(m)}[u]
=\sum_{k\geq 0}(-1)^{k}\binom{n}{k}[a]_{(n-k)}([c]_{(m+k)}[u])
\end{eqnarray}
\end{lemma}
\noindent{\bf Proof.}
(i)\  We prove that identity (\ref{e2.3}) is true when $[a]=b,\ |c|=1$.

Let $[c]=D^{i}b'$. If $0\leq n<N(b,b')$, identity (\ref{e2.3}) follows from the definition.

Let $n\geq N(b,b')$. When $i=0$, the left hand side of identity (\ref{e2.3}) is equal to $0$,
while the right hand side of identity (\ref{e2.3}) contains the summand
$$
b_{(n)}(b'_{(m)}[u])
=-\sum_{k\geq 1}(-1)^{s}\binom{n}{k}b_{(n-k)}(b'_{(m+k)}[u]).
$$
Hence the right hand side of identity (\ref{e2.3}) is $0$ as well.

Induction on $i$. Suppose that $i\geq 1$. The left hand side of identity (\ref{e2.3}) is equal to
\begin{align*}
(b_{(n)}D^{i}b')_{(m)}[u]
&=D(b_{(n)}D^{i-1}b')_{(m)}[u]+n(b_{(n-1)}D^{i-1}b')_{(m)}[u] \\
&=-m(b_{(n)}D^{i-1}b')_{(m-1)}[u]+n(b_{(n-1)}D^{i-1}b')_{(m)}[u].
\end{align*}
The right hand side of identity (\ref{e2.3}) is equal to
\begin{align*}
&\sum_{k\geq 0}(-1)^{k}\binom{n}{k}b_{(n-k)}(D^{i}b'_{(m+k)}[u])
=-\sum_{k\geq 0}(-1)^{k}(m+k)\binom{n}{k}b_{(n-k)}(D^{i-1}b'_{(m+k-1)}[u]) \\
=&-m\sum_{k\geq 0}(-1)^{k}\binom{n}{k}b_{(n-k)}(D^{i-1}b'_{(m-1+k)}[u])+
\sum_{k\geq 1}(-1)^{k+1}k\binom{n}{k}b_{(n-k)}(D^{i-1}b'_{(m+k-1)}[u]) \\
=&-m(b_{(n)}D^{i-1}b')_{(m-1)}[u]+\sum_{k\geq 0}(-1)^{k}(k+1)\binom{n}{k+1}b_{(n-k-1)}(D^{i-1}b'_{(m+k)}[u]) \\
=&-m(b_{(n)}D^{i-1}b')_{(m-1)}[u]+n\sum_{k\geq 0}(-1)^{k}\binom{n-1}{k}b_{(n-1-k)}(D^{i-1}b'_{(m+k)}[u])\\
=&-m(b_{(n)}D^{i-1}b')_{(m-1)}[u]+n(b_{(n-1)}D^{i-1}b')_{(m)}[u].
\end{align*}
Hence, identity (\ref{e2.3}) is true for $[a]=b,\ |c|=1$.

(ii)\  We  prove that identity (\ref{e2.3}) is true for any $c\in T,\ [a]=b$.

Induction on $|c|$. We have showed the result for $|c|=1$. Suppose $|c|>1$ and $[c]=b_{1(p)}[c_1]$.
We only need to consider $n\geq N(b,b_1)$. Following the definition, we have
$$
b_{(n)}(b_{1(p)}[c_1])=-\sum_{k\geq 1}(-1)^{k}\binom{n}{k}b_{(n-k)}(b_{1(p+k)}[c_1]).
$$
Now induction on $n$. Let $n=N(b,b_1)=:N$. Then
\begin{align*}
A:&=(b_{(N)}[c])_{(m)}[u]=-\sum_{k\geq 1}(-1)^{k}\binom{N}{k}(b_{(N-k)}(b_{1(p+k)}[c_1]))_{(m)}[u] \\
&=-\sum\limits_{\scriptstyle k\geq 1,\ t\geq 0}(-1)^{k+t}\binom{N}{k}\binom{N-k}{t}b_{(N-k-t)}((b_{1(p+k)}[c_1])_{(m+t)}[u])\\
&=\sum_{\scriptstyle k\geq 1,\ t,r\geq 0}(-1)^{k+t+r+1}\binom{N}{k}\binom{N-k}{t}\binom{p+k}{r}b_{(N-k-t)}(b_{1(p+k-r)}([c_1]_{(m+t+r)}[u])).
\end{align*}

We denote the right hand side of identity (\ref{e2.3}) is $A_1+A_2$ where
$A_1=b_{(N)}((b_{1(p)}[c_1])_{(m)}[u])$,\ $A_2=\sum_{k\geq 1}(-1)^{k}\binom{N}{k}b_{(N-k)}((b_{1(p)}[c_1])_{(m+k)}[u])$.
Then
\begin{eqnarray*}
A_1&=&\sum_{t\geq 0}(-1)^{t}\binom{p}{t}b_{(N)}(b_{1(p-t)}([c_1]_{(m+t)}[u])) \\
&=&-\sum\limits_{t\geq 0,r\geq 1}(-1)^{t+r}\binom{p}{t}\binom{N}{r}b_{(N-r)}(b_{1(p-t+r)}([c_1]_{(m+t)}[u])),
\\
A_2&=&\sum_{k\geq 1,t\geq 0}(-1)^{k+t}\binom{N}{k}\binom{p}{t}b_{(N-k)}(b_{1(p-t)}([c_1]_{(m+k+t)}[u])).
\end{eqnarray*}

We make a transformation
$$
(i,\ j,\ l)=(N-k-t,\ p+k-r,\ m+t+r),
$$
where $i,\ j,\ l$ are nonnegative integers, and so $i+j+l=N+m+p$.
Then $A$ becomes a sum of the expressions
\begin{eqnarray*}
A&=&\sum_{i,j,l\geq 0,\ k\geq 1}(-1)^{k+l+m+1}\binom{N}{k}\binom{N-k}{N-k-i}\binom{p+k}{p+k-j}b_{(i)}(b_{1(j)}([c_1]_{(l)}[u]))\\
&=&\sum_{i,j,l\geq 0,\ k\geq 1}(-1)^{k+l-m+1}\binom{N}{k}\binom{N-k}{i}\binom{p+k}{j}b_{(i)}(b_{1(j)}([c_1]_{(l)}[u])) \\
&=&\sum_{i,j,l\geq 0,\ k\geq 1}(-1)^{k+N+p-i-j+1}\binom{N}{i}\binom{N-i}{k}\binom{p+k}{j}b_{(i)}(b_{1(j)}([c_1]_{(l)}[u])) \\
&=&\sum_{i,j,l\geq 0}(-1)^{p+j+1}\binom{N}{i}(\sum_{k\geq 1}(-1)^{N-i-k}\binom{N-i}{k}\binom{p+k}{j})b_{(i)}(b_{1(j)}([c_1]_{(l)}[u]))\\
&=&\sum_{i,j,l\geq 0}(-1)^{p+j+1}\binom{N}{i}(\binom{p}{j-N+i}-(-i)^{N-i}\binom{p}{j})b_{(i)}(b_{1(j)}([c_1]_{(l)}[u]))\\
&=&\sum_{i,j,l\geq 0}\binom{N}{i}((-1)^{p+j+1}\binom{p}{l-m}+(-1)^{l-m}\binom{p}{j})b_{(i)}(b_{1(j)}([c_1]_{(l)}[u])).
\end{eqnarray*}

Next, do a similar transformation
$$
(i,\ j,\ l)=(N-r,\ p+r-t,\ m+t),
$$
where $i,\ j,\ l$ are nonnegative integers, and so $i+j+l=N+m+p$.
Then
\begin{align*}
A_1&=\sum_{i,j,l\geq 0}(-1)^{p+j+1}\binom{N}{N-i}\binom{p}{l-m}b_{(i)}(b_{1(j)}([c_1]_{(l)}[u]))\\
&=\sum_{i,j,l\geq 0}(-1)^{p+j+1}\binom{N}{i}\binom{p}{l-m}b_{(i)}(b_{1(j)}([c_1]_{(l)}[u])).
\end{align*}
Do another transformation
$$
(i,\ j,\ l)=(N-s,\ p-t,\ m+s+t),
$$
where $i,\ j,\ l$ are nonnegative integers, and so $i+j+l=N+m+p$.
Then
\begin{align*}
A_2&=\sum_{i,j,l\geq 0}(-1)^{l+m}\binom{N}{N-i}\binom{p}{p-j}b_{(i)}(b_{1(j)}([c_1]_{(l)}[u]))\\
&=\sum_{i,j,l\geq 0}(-1)^{l-m}\binom{N}{i}\binom{p}{j}b_{(i)}(b_{1(j)}([c_1]_{(l)}[u])).
\end{align*}
Thus $A=A_1+A_2$. Assume that $n>N(b,b_1)$, we just repeat the argument of $n=N(b,b_1)$.
So identity (\ref{e2.3}) is true for $[a]=b$.

(iii)\ We  prove that identity (\ref{e2.3}) holds for $[a]=D^{j}b,\ j\geq 1$.

By definition, the left hand side of identity (\ref{e2.3}) is equal to
$$
(D^{j}b_{(n)}[c])_{(m)}[u]=(-1)^{j}\frac{n!}{(n-j)!}(b_{(n-j)}[c])_{(m)}[u].
$$
The right hand side of identity (\ref{e2.3}) is equal to
\begin{align*}
& \ \ \ \ \sum_{k\geq 0}(-1)^{k}\binom{n}{k}D^{j}b_{(n-k)}([c]_{(m+k)}[u]) \\
&=-\sum_{k\geq 0}(-1)^{k+j}\binom{n}{k}\frac{(n-k)!}{(n-k-j)!}b_{(n-k-j)}([c]_{(m+k)}[u]) \\
&=(-1)^{j}\frac{n!}{(n-j)!}\sum_{k\geq 0}(-1)^{k}\binom{n-j}{k}b_{(n-j-k)}([c]_{(m+k)}[u]).
\end{align*}
Therefore, identity (\ref{e2.3}) holds for $|a|=1$.

(iv)\ Assume that $|a|>1$. In this case, we write $[a]=b_{(q)}[a_1]$.
Then the left hand side of identity (\ref{e2.3}) is equal to
\begin{align*}
&\ \ \ \ ((b_{(q)}[a_1])_{(n)}[c])_{(m)}[u]=\sum\limits_{k\geq 0}(-1)^{k}\binom{q}{k}(b_{(q-k)}([a_1]_{(n+k)}[c]))_{(m)}[u] \\
&=\sum\limits_{k,t\geq 0}(-1)^{k+t}\binom{q}{k}\binom{q-k}{t}b_{(q-k-t)}(([a_1]_{(n+k)}[c])_{(m+t)}[u]) \\
&=\sum\limits_{k,t,r\geq 0}(-1)^{k+t+r}\binom{q}{k}\binom{q-k}{t}\binom{n+k}{r}b_{(q-k-t)}([a_1]_{(n+k-r)}([c]_{(m+t+r)}[u]))=:A.
\end{align*}
Do a transformation on the indices
$$
(i,\ j,\ l)=(q-k-t,\ n+k-r,\ m+t+r),
$$
where $i,\ j,\ l$ are nonnegative integers, and so $i+j+l=n+m+q$.
Then
\begin{align*}
A&=\sum_{i,j,l,k\geq 0}(-1)^{k+l+m}\binom{q}{k}\binom{q-k}{q-k-i}\binom{n+k}{n+k-j}b_{(i)}([a_1]_{(j)}([c]_{(l)}[u])) \\
&=\sum_{i,j,l,k\geq 0}(-1)^{k+n+q-i-j}\binom{q}{k}\binom{q-k}{i}\binom{n+k}{j}b_{(i)}([a_1]_{(j)}([c]_{(l)}[u])) \\
&=\sum_{i,j,l\geq 0}(-1)^{n-j}\binom{q}{i}(\sum_{k\geq 0}(-1)^{q-i-k}\binom{q-i}{k}\binom{n+k}{j})b_{(i)}([a_1]_{(j)}([c]_{(l)}[u])) \\
&=\sum_{i,j,l\geq 0}(-1)^{n-j}\binom{q}{i}\binom{n}{j-q+i}b_{(i)}([a_1]_{(j)}([c]_{(l)}[u]))\\
&=\sum_{i,j,l\geq 0}(-1)^{n-j}\binom{q}{i}\binom{n}{n+m-l}b_{(i)}([a_1]_{(j)}([c]_{(l)}[u]))\\
&=\sum_{i,j,l\geq 0}(-1)^{n-j}\binom{q}{i}\binom{n}{l-m}b_{(i)}([a_1]_{(j)}([c]_{(l)}[u])).
\end{align*}

Denote the right hand side of identity (\ref{e2.3}) as $A'$. Then
\begin{align*}
A'&=\sum\limits_{k\geq 0}(-1)^{k}\binom{n}{k}(b_{(q)}[a_1])_{(n-k)}([c]_{(m+k)}[u]) \\
&=\sum\limits_{k,t\geq 0}(-1)^{k+t}\binom{n}{k}\binom{q}{t}b_{(q-t)}([a_1]_{(n-k+t)}([c]_{(m+k)}[u]))\\
&=\sum\limits_{i,j,l\geq 0}(-1)^{n-j}\binom{n}{l-m}\binom{q}{q-i}b_{(i)}([a_1]_{(j)}([c]_{(l)}[u]))=A,
\end{align*}
where
$
(i,\ j,\ l)=(q-t,\ n-k+t,\ m+k).
$
Therefore, we complete the proof of identity (\ref{e2.3}).\ \ \ \ $\square$

\begin{theorem}\label{th3.1}
Let the notation be as above. Then $span_{\mathbf{k}}[U]$ is the free  module generated by $Y$ over the free associative
conformal algebra $C(B,N)$.

We denote $span_{\mathbf{k}}[U]$ by $mod_{C(B,N)}(Y)$, the double free associative conformal module.
\end{theorem}
\noindent{\bf Proof.}
By Lemmas \ref{l2.1}--\ref{l2.5}, $span_{\mathbf{k}}[U]$ is a  module over free associative
conformal algebra $C(B,N)$.

For any $C(B,N)$-module $M$ and any map $\varepsilon$ satisfying $[a]_{(n)}\varepsilon(y)=0$ for all $a\in T, y\in Y, n\geq N([a],y)$,
we define the $C(B,N)$-module homomorphism
$$
\varphi: \ mod_{C(B,N)}(Y)\longrightarrow \ _{C(B,N)}M, \ \ \
[u]\longmapsto \ [u]|_{y\mapsto\varepsilon(y)}
$$
such that the following diagram is commutative:

\setlength{\unitlength}{1cm}
\begin{picture}(7, 3)
\put(4.2,2.3){\vector(1,0){1.7}} \put(4.1, 2.0){\vector(0,-1){1.3}}
\put(5.9,2.1){\vector(-1,-1){1.6}} \put(3.9,0.2){$M$}
\put(6,2.2){$mod_{C(B,N)}(Y)$} \put(3.9,2.2){$Y$} \put(4.9, 2.4){$i$}
\put(3.7,1.3){$\varepsilon$} \put(5.3,1){$\exists!\varphi$}
\end{picture}

Thus, the result is true.\ \ \ \ $\square$

\subsection{Free associative conformal $C$-modules}

Let $C=(C, N(-,-), (n), n\in\mathbb{Z}_{\geq 0})$ be an arbitrary associative conformal algebra.
Then $C$ is an epimorphic image of some free associative conformal algebra $C(B, N)$. Thus, $C$ has an expression
$$
C=C(B,N|S):=C(B,\ N)/Id(S)
$$
generated by $B$ with defining relations $S$, where $Id(S)$ is the ideal of
$C(B,N)$ generated by $S$.

Let $S\subset C(B,N),\ C=C(B,N|S)$ and $mod_{C(B,N)}(Y)$ be the  $C(B,N)$-module constructed as above. Denote
\begin{eqnarray*}
R:&=&\{s_{(m)}[u]\mid  s\in S, u\in U, m\geq 0\},\\
D^{\omega}(Y):&=&\{D^{i}y\ | \ i\geq 0,\ y\in Y\}.
\end{eqnarray*}
Then $subm(R)=\sum\limits_{m\geq 0} Id(S)_{(m)}D^{\omega}(Y)$. Thus,
$mod_{C(B,N)}(Y|R)$ is also a $C$-module if we define: for any
$f+Id(S)\in C,\ n\geq 0,\ h+subm(R)\in mod_{C(B,N)}(Y|R)$,
$$
(f+Id(S))_{(n)}(h+subm(R)):=f_{(n)}h+subm(R).
$$

For any left $C$-module $_CM$, we also can regard $_CM$ as a $C(B,N)$-module
in a natural way: for any $f\in C(B,N), \ n\geq 0,\ v\in M$,
$$
f_{(n)}v:=(f+Id(S))_{(n)}v.
$$

\begin{theorem}\label{pro1}
Let $C=C(B,N|S)$ be an arbitrary associative conformal algebra, $Y$ a set and
$
R=\{s_{(m)}[u]\mid s\in S, u\in U, m\geq 0\}.
$
Then $mod_{C(B,N)}(Y|R)$ is a free $C$-module generated by $Y$.

We denote $mod_{C(B,N)}(Y|R)$ by $mod_{C(B,N|S)}(Y)$.
\end{theorem}
\noindent{\bf Proof.} \
Let $_{C}M$ be a $C$-module and $\varepsilon: Y\rightarrow M$ be a map such that $c_{(n)}\varepsilon(y)=0$
for any $n\geq N(c,y),\ c\in C,\ y\in Y$. By Theorem \ref{th3.1},
there exists a unique $C(B,N)$-module homomorphism $\phi: mod_{C(B,N)}(Y)\rightarrow _C\!\!M$,
such that $\phi i=\varepsilon$. Let $\pi$ be the natural $C(B,N)$-module homomorphism from $mod_{C(B,N)}(Y)$ to $mod_{C(B,N)}(Y|R)$.
Obviously, for any $s_{(m)}[u]\in R$,
$$
\phi(s_{(m)}[u])=s_{(m)}\phi([u])=(s+Id(S))_{(m)}\phi([u])=0.
$$
Hence, there exists a unique $C(B,N)$-module homomorphism $\varphi:   mod_{C(B,N)}(Y|R)\rightarrow _C\!\!M$
satisfying $\varphi\pi=\phi$.

\setlength{\unitlength}{1cm}
\begin{picture}(7, 3)
\put(4.2,2.3){\vector(1,0){1.7}} \put(4.1, 2.0){\vector(0,-1){1.3}}
\put(5.9,2.1){\vector(-1,-1){1.6}} \put(11,2.0){\vector(-4,-1){6.6}}
\put(8.7,2.3){\vector(1,0){2.1}}\put(3.5,0.2){$_CM$}
\put(6,2.2){$mod_{C(B,N)}(Y)$} \put(3.9,2.2){$Y$} \put(11,2.2){$mod_{C(B,N)}(Y|R)$}
\put(4.9, 2.4){$i$}\put(3.7,1.3){$\varepsilon$} \put(9.7,2.4){$\pi$}
\put(5.3,1){$\exists!\phi$}\put(9,1){$\exists!\varphi$}
\end{picture}

Noting that $\varphi$ is also a $C$-module homomorphism, we complete the proof.   \ \ \ \  $\square$

\ \

Let $_CM$ be a $C$-module. Then $_CM$ has an expression $_CM=mod_{C(B,N|S)}(Y|Q)$.

By noting that $ mod_{C(B,N|S)}(Y)=mod_{C(B,N)}(Y|R)$, we may assume that $Q\subseteq mod_{C(B,N)}(Y|R)$.
Let $Q_{1}=\{f\in mod_{C(B,N)}(Y)\mid f+subm(R)\in Q \}$. Then we have

\begin{proposition}\label{p3}
Let the notation as above. Then
$$
mod_{C(B,N|S)}(Y|Q)\cong mod_{C(B,N)}(Y|R\cup Q_{1})
$$
as $C(B,N)$-modules and  as $C(B,N|S)$-modules.
\end{proposition}

\section{Composition-Diamond lemma for associative conformal modules}

We give Composition-Diamond lemma for double free associative
conformal module $mod_{C(B,N)}(Y)$ step by step.

In this section, we fix the free $C(B,N)$-module $mod_{C(B,N)}(Y)$ generated by $Y$ with a uniform bounded locality function $N$.
Recall that
$$
T=\{b_{1(n_{1})}\cdots b_{k(n_{k})}D^{i}b_{k+1}\mid  b_{j},b_{k+1}\in B,
0\leq n_{j}<N, 1\leq j\leq k, i, k\geq 0\},
$$
$$
U=\{b_{1(n_{1})}b_{2(n_{2})}\cdots b_{k(n_{k})}D^{i}y\mid
b_{j}\in B,\ y\in Y,\ 0\leq n_{j}<N,\ 1\leq j\leq k,\ i, k\geq 0 \}
$$
and $[U]=\{[u]\mid u\in U\}$ is a $\mathbf{k}$-basis of $mod_{C(B,N)}(Y)$, where  $[\ ]$ means right
normed bracketing.

For $j\geq 0,\ u\in U$, $uD^{j}$ means that we apply
$D^{j}$ to the last letter of $u$.

A word in $mod_{C(B,N)}(Y)$ is a polynomial of the form
$(u)$ (with some bracketing of $u$), where
$u= b_{1(n_{1})}\cdots b_{k(n_{k})}D^{i}y,\ b_{j}\in B,\ y\in Y,\ n_{j}, i, k\geq 0,\ 1\leq j\leq k $.
Since $[U]$ is  a $\mathbf{k}$-basis of $mod_{C(B,N)}(Y)$, each word $(u)$ in $mod_{C(B,N)}(Y)$ is a linear combination of some
elements in $[U]$.

\subsection{A monomial ordering}

Let $Y,B$ be well-ordered sets. We order elements of $U$ according
to the lexicographical ordering of their weights. For any
$u=b_{1(n_{1})}b_{2(n_{2})}\cdots
b_{k(n_{k})}D^{i}y\in U$, denote by
$$
wt(u):=(|u|,b_{1},n_{1},\cdots, b_{k},n_{k},y,i),
$$
where $|u|=k+1$ is the length of $u$.
Then for any $u, v\in U$, define
$$
u>v\Leftrightarrow wt(u)>wt(v)\ \ \ \ \mbox{ lexicographically.}
$$
It is clear that such an ordering is a well ordering on $U$. We will use this ordering in the sequel.

For $f\in mod_{C(B,N)}(Y)$, the leading term of $f$ is denoted by $\bar f$ and $\bar f\in U$.
So $f=\alpha_{\bar f}[\bar f]+\sum_{i}\alpha_{i}[u_i]$, $u_i<\bar f$.
We will call $f$ monic if $\alpha_{\bar f}=1$.

\begin{lemma}\label{l3.2}
Suppose that $a\in T$ is $D$-free, $0\leq n<N$, $u, v\in U$. Then
\begin{enumerate}
 \item[(i)]\ $\overline{[a]_{(n)}[u]}=a_{(n)}u$;
 \item[(ii)]\ $u<v\Longrightarrow \overline{[a]_{(n)}[u]}<\overline{[a]_{(n)}[v]}$;
 \item[(iii)]\ $\overline{D^{i}([u])}=uD^{i}$, $i\in \mathbb{Z}_{\geq 0}$;
 \item[(iv)]\ $u<v\Longrightarrow \overline{D^{i}([u])}<
 \overline{D^{i}([v])}$, $ i\in \mathbb{Z}_{\geq 0}$.
\end{enumerate}
\end{lemma}
\noindent{\bf Proof.} (i) \  Induction on $|a|$. If $|a|=1$,
then $[a]=b$ for some $b\in B$ and
$$
\overline{[a]_{(n)}[u]}=\overline{b_{(n)}[u]}=b_{(n)}u.
$$
Assume that $|a|>1$ and $[a]=b_{(m)}[a_{1}]$, where
$a_{1}\in T$ is $D$-free. Then
$$
[a]_{(n)}[u]=(b_{(m)}[a_{1}])_{(n)}[u]=\sum_{k=0}^{m}
(-1)^{k}\binom{m}{k}b_{(m-k)}([a_{1}]_{(n+k)}[u]).
$$
Since $m<N$,
\begin{align*}
\overline{[a]_{(n)}[u]}
=\overline{b_{(m)}([a_{1}]_{(n)}[u])}=b_{(m)}(\overline{[a_{1}]_{(n)}[u]})
=b_{(m)}a_1{_{(n)}}u=a_{(n)}u.
\end{align*}

(ii) \  This part follows from (i).

(iii) \  Induction on $|u|$. The result is obvious for
$|u|=1$. Let $|u|>1$, $[u]=b_{(m)}[u_{1}]$. Then
$D^{i}([u])=\sum\limits_{k\geq 0}(-1)^{k}\binom{i}{k}\frac{m!}{(m-k)!}b_{(m-k)}D^{i-k}([u_{1}])$.
Therefore,
$$
\overline{D^{i}([u])}
=\overline{b_{(m)}D^{i}([u_{1}])}=b_{(m)}\overline{D^{i}([u_{1}])}
=b_{(m)}u_{1}D^{i}
=uD^{i}.
$$

(iv) \   This part follows from (iii).   \ \ \ \  $\square$

\ \

Thus, the ordering $>$ on $U$ is a monomial ordering in a sense of
(ii) and (iv) in Lemma \ref{l3.2}.

\subsection{$S$-words and normal $S$-words}

Let $S\subset mod_{C(B,N)}(Y)$ be a set of
monic polynomials, $(u)$ be a word in $mod_{C(B,N)}(Y)$, where
$u= b_{1(n_{1})}\cdots b_{k(n_{k})}D^{i}y,\ b_{j}\in B,\ y\in Y,\ n_{j}, i, k\geq 0,\ 1\leq j\leq k $.
We define $S$-word $(u)_{D^{i}s}$ for any $s\in S$ by induction.

(i)\ $(D^{i}s)_{D^{i}s}=D^{i}s$ is an $S$-word of $S$-length 1;

(ii)\ If $(u)_{D^{i}s}$ is an $S$-word of $S$-length
$k$, and $(a)$ is any word in $C(B,N)$ of length $l$,
then $(a)_{(n)}(u)_{D^{i}s}$ is an $S$-words of $S$-length $k+l$.

The $S$-length of an $S$-word $(u)_{D^{i}s}$ will be
denoted by $|u|_{D^{i}s}$.

An $S$-word $\lfloor u\rfloor_{D^{i}s}:=b_{1(n_{1})}(b_{2(n_{2})}(\cdots (b_{k(n_{k})}D^{i}s)\cdots ))$
 is called a right normed $S$-word, and it is a normal $S$-word, denoted by $[u]_{D^{i}s}$, if each $n_j<N$.

\begin{lemma}
Let $[u]_{D^{i}s}$ be a normal $S$-word. Then $\overline
{[u]_{D^{i}s}}=u|_{_{D^{i}s\mapsto\overline{s}D^{i}}}$.
\end{lemma}
\noindent{\bf Proof.}  Let $s=[\overline{s}]+\sum_{j}\alpha_{j}
[v_{j}]$, where $v_{j}<\overline{s}$.
Then $D^{i}s=D^{i}([\overline{s}])+\sum_{j}\alpha_{j}D^{i}([v_{j}])$,
and  applying Lemma \ref{l3.2}, we get $\overline{D^{i}([v_{j}])}<\overline{D^{i}([\overline{s}])}$.
Hence $\overline{D^{i}s}=\overline{D^{i}([\overline{s}])}=
\overline{s}D^{i}$.

If $|u|_{D^{i}s}=1$, then $[u]_{D^{i}s}=D^{i}s$ and we are done.
So assume that $|u|_{D^{i}s}>1$. Then $[u]_{D^{i}s}=b_{(n)}[v]_{D^{i}s}$,
so
\begin{align*}
\overline{[u]_{_{D^{i}s}}}
=\overline{b_{(n)}[v]_{_{D^{i}s}}}
=b_{(n)}\overline{[v]_{_{D^{i}s}}}
=b_{(n)}v|_{_{D^{i}s\mapsto \overline{s}D^{i}}}
=u|_{_{D^{i}s\mapsto \overline{s}D^{i}}}.
\end{align*}
We complete the proof by induction on $|u|_{D^{i}s}$. \ \ \ \  $\square$

\begin{lemma}\label{l3.4}
Any $S$-word $(u)_{_{D^{i}s}}$ can be presented as a linear combination
of right normed $S$-word such that the length of the leading term of
each term is less or equal to  $|\overline{(u)_{_{D^{i}s}}}|$.
\end{lemma}
\noindent{\bf Proof.}  Induction on $|u|_{_{D^{i}s}}$.
The result holds trivially if $|u|_{_{D^{i}s}}=1$.
Let $|u|_{_{D^{i}s}}>1$. We assume that
\begin{eqnarray}\label{e3.5}
(u)_{_{D^{i}s}}=(a)_{(m)}\lfloor v\rfloor_{_{D^{i}s}}
\end{eqnarray}
where $\lfloor v\rfloor_{_{D^{i}s}}$ is a right normed $S$-word.
Next, induction on $|a|$. If $|a|=1$, then $(a)=D^{j}b$
for some $b\in B$ and $j\geq 0$, so
$$
(a)_{(m)}\lfloor v\rfloor_{_{D^{i}s}}=D^{j}b_{(m)}\lfloor v\rfloor_{_{D^{i}s}}=\left\{
\begin{aligned}
&0, \ \ \ \ \ \ \ \ \ \ \ \ \ \ \ \ \ \ \ \ \ \ \ \ \ \ \ \ \ \ \ \ \ \ \ \   j> m , \\
&(-1)^{j}\frac{m!}{(m-j)!}b_{(m-j)}\lfloor v\rfloor_{_{D^{i}s}}, \ \ j\leq m.
\end{aligned}
\right.
$$
Thus, (\ref{e3.5}) is right normed and we are done. Let $|a|>1$ and
$(a)=(a_{1})_{(n)}(a_{2})$. Then
$((a_{1})_{(n)}(a_{2}))_{(m)}\lfloor v\rfloor_{_{D^{i}s}}=\sum\limits_{k=0}^{n}
(-1)^{k}\binom{n}{k}(a_{1})_{(n-k)}((a_{2})_{(m+k)}\lfloor v\rfloor_{_{D^{i}s}})$.
Now, the result follows from the induction on $|a|$.  \ \ \ \
$\square$

\subsection{Compositions}

Let $S\subset mod_{C(B,N)}(Y)$ with each polynomial in $S$ monic, $w\in U$ and $f,g\in S$.

We have three kinds of compositions.

$\bullet$ If $w=\bar f=\overline{[u]_{D^{i}g}}$ and $i\geq 0$, then define
$$
(f,g)_w=f-[u]_{D^{i}g},
$$
which is a composition of inclusion.

$\bullet$ If $w=\bar fD^{i}=a_{(n)}\bar g, i>0$ and $a\in T$ with $a\ D$-free, then define
$$
(f,g)_w=D^{i}f-[a_{(n)}g],
$$
which is a composition of intersection.

$\bullet$ If $b\in B$ and $n\geq N$, then  $b_{(n)}f$ is
referred to as a composition of left multiplication.

\ \

Let $S\subset mod_{C(B,N)}(Y)$ be a set of monic polynomials and $h\in mod_{C(B,N)}(Y)$.
Then $h$ is said to be trivial modulo $S$, denoted by

\begin{eqnarray}\label{e3.3}
h\equiv0\ mod(S)
\end{eqnarray}
if $h=\sum_{i}\alpha_{i}[u_{i}]_{_{D^{l_i}s_{i}}}$, where each $[u_{i}]_{_{D^{l_i}s_{i}}}$ is a normal $S$-word and
$\overline{[u_{i}]_{_{D^{l_i}s_{i}}}}\leq \bar{h}$. For $h_1,h_2\in mod_{C(B,N)}(Y)$,
$h_1\equiv h_2\ mod(S)$ means $h_1-h_2\equiv 0\ mod(S)$.

The  set $S$ is
called a Gr\"{o}bner-Shirshov basis in $mod_{C(B,N)}(Y)$ if all compositions of elements
of $S$ are trivial modulo $S$. In particular, if $S$ is $D$-free, then we call $S$ a $D$-free Gr\"{o}bner-Shirshov basis.

The set $S$ is said to be closed under the composition of left multiplication if any composition of left multiplication of $S$ is trivial modulo $S$.
That $S$ is closed under the composition of inclusion
and intersection is similarly defined.

\begin{lemma}\label{l3.3}
Let $S\subset mod_{C(B,N)}(Y)$ with each polynomial in $S$ monic and $[u]_{D^{i}s}$ a normal $S$-word. Then
\begin{enumerate}
\item[(i)]\ $D^{j}([u]_{D^{i}s})=[u]_{D^{i+j}s}+\sum_{t}\beta_{t}[w_{t}]_{_{D^{l_t}s_{t}}}$,
where each $[w_{t}]_{_{D^{l_t}s_{t}}}$ is a normal $S$-word and $\overline{[w_{t}]_{_{D^{l_t}s_{t}}}}<\overline{[u]_{D^{i+j}s}}$.
\item[(ii)]\ If $S$ is closed under the composition of left multiplication,
then for any $b\in B,\ n\geq N$, we have
$b_{(n)}[u]_{D^{i}s}\equiv 0\ mod(S)$.
\end{enumerate}
\end{lemma}
\noindent{\bf Proof.} (i) \ Induction on $|u|_{D^{i}s}$. Let $|u|_{D^{i}s}=1$ and $[u]_{D^{i}s}=D^{i}s$.
Then
$$
D^{j}([u]_{D^{i}s})=D^{j}(D^{i}s)=D^{i+j}s=[u]_{D^{i+j}s}.
$$
Assume that $|u|_{D^{i}s}>1$ and $[u]_{D^{i}s}=b_{(n)}[v]_{D^{i}s}$. So
\begin{align*}
D^{j}([u]_{D^{i}s})
&=D^{j}(b_{(n)}[v]_{D^{i}s})=\sum_{p\geq 0}(-1)^{p}\binom{j}{p}
\frac{n!}{(n-p)!}b_{(n-p)}D^{j-p}([v]_{D^{i}s})  \\
&=b_{(n)}D^{j}([v]_{D^{i}s})+\sum_{p\geq 1}(-1)^{p}\binom{j}{p}
\frac{n!}{(n-p)!}b_{(n-p)}D^{j-p}([v]_{D^{i}s})\\
&=b_{(n)}[v]_{D^{i+j}s}+\sum_{t}\beta_{t}b_{(n)}[w_{t}]_{_{D^{l_t}s_{t}}}
+\sum_{p\geq 1}\alpha_{p}b_{(n-p)}([v]_{D^{i+j-p}s}+\sum_{t'}\gamma_{t'}[w'_{t'}]_{_{D^{l_{t'}}s_{t'}}})\\
&=[u]_{D^{i+j}s}+\sum_{t}\beta_{t}b_{(n)}[w_{t}]_{_{D^{l_t}s_{t}}}
+\sum_{p\geq 1}\alpha_{p}b_{(n-p)}([v]_{D^{i+j-p}s}+\sum_{t'}\gamma_{t'}[w'_{t'}]_{_{D^{l_{t'}}s_{t'}}}),
\end{align*}
where $\overline{[w_{t}]_{_{D^{l_t}s_{t}}}}< \overline{[v]_{D^{i+j}s}}$,\ $\overline{[w'_{t'}]_{_{D^{l_{t'}}s_{t'}}}}< \overline{[v]_{D^{i+j-p}s}}$,\
$\alpha_p=(-1)^{p}\binom{j}{p}\frac{n!}{(n-p)!}$.
Thus $\overline{b_{(n)}[w_{t}]_{_{D^{l_t}s_{t}}}}<\overline{b_{(n)}[v]_{D^{i+j}s}}=\overline{[u]_{D^{i+j}s}}$
and $\overline{b_{(n-p)}D^{j-p}([v]_{D^{i}s})}=b_{(n-p)}\overline{[v]_{D^{i+j-p}s}}<b_{(n)}\overline{[v]_{D^{i+j}s}}
=\overline{[u]_{D^{i+j}s}}$ for $p\geq 1$.
By induction, we have
$D^{j}([u]_{D^{i}s})\equiv[u]_{D^{i+j}s}\ mod(S,w)$.

(ii) \ Let $|u|_{D^{i}s}=1$ and $n=N$. Then
\begin{align*}
b_{(N)}D^{i}s
&=D^{i}(b_{(N)}s)-\sum_{t\geq 1}(-1)^{t}
\binom{i}{t}\frac{N!}{(N-t)!}b_{(N-t)}D^{i-t}s \\
&=D^{i}(b_{(N)}s)-\sum_{t\geq 1}\alpha_{t}[u_{t}]_{D^{i-t}s},
\end{align*}
where $\alpha_{t}=(-1)^{t}\binom{i}{t}\frac{N!}{(N-t)!}$,
and $[u_{t}]_{D^{i-t}s}=b_{(N-t)}D^{i-t}s$.
We may assume that $b_{(N)}s=\sum\limits_{j}\alpha_{j}[v_{j}]_{D^{l_{j}}s_{j}}$,
where $\overline{[v_{j}]_{D^{l_{j}}s_{j}}}\leq\overline{b_{(N)}s}$.
By (i), we get
\begin{align*}
D^{i}(b_{(N)}s)
&=\sum_{j}\alpha_{j}D^{i}([v_{j}]_{D^{l_{j}}s_{j}})=
\sum_{j}\alpha_{j}([v_{j}]_{D^{l_{j}+i}s_{j}}
+\sum_{q}\beta_{j_{q}}[w_{j_{q}}]_{D^{l_{j_{q}}}s_{j_{q}}}) \\
&=\sum_{j}\alpha_{j}[v_{j}]_{D^{l_{j}+i}s_{j}}+\sum_{j,j_{q}}
\alpha_{j}\beta_{j_{q}}[w_{j_{q}}]_{D^{l_{j_{q}}}s_{j_{q}}},
\end{align*}
where $\overline{[w_{j_{q}}]_{D^{l_{j_{q}}}s_{j_{q}}}}<\overline{[v_{j}]_{D^{l_{j}+i}s_{j}}} $
and $\overline{[v_{j}]_{D^{l_{j}+i}s_{j}}}
=\overline{D^{i}([v_{j}]_{D^{l_{j}}s_{j}})}\leq \overline{D^{i}(b_{(N)}s)}$.
By induction on $n$, the result
is true for $|u|_{D^{i}s}=1$.

Assume that $|u|_{D^{i}s}>1$ and the result holds for any $n\geq N$
when the $S$-length of a normal $S$-word is less than
$|u|_{D^{i}s}$. Let $[u]_{D^{i}s}=b_{1(n_{1})}[v]_{D^{i}s}$. Then
$$
b_{(n)}[u]_{D^{i}s}=b_{(n)}(b_{1(n_{1})}[v]_{D^{i}s})=
-\sum_{t\geq 1}(-1)^{t}\binom{n}{t}b_{(n-t)}(b_{1(n_{1}+t)}[v]_{D^{i}s}).
$$
Hence, the result follows from the induction on $n$.
\ \ \ \  $\square$

\begin{lemma}\label{ll3.4}
Let $S$ be a subset of monic polynomials, closed under the composition of
left multiplication.
Then
any $S$-word $(u)_{D^{i}s}$ has a presentation
$
(u)_{D^{i}s}=\sum\limits_{j}\alpha_{j}[u_j]_{D^{l_j}s_j} ,
$
where $[u_j]_{D^{l_j}s_j}$ is a normal $S$-word for each $j$.
\end{lemma}
\noindent{\bf Proof.}  Due to Lemma \ref{l3.4}, we may assume
that $(u)_{D^{i}s}$ is right normed, i.e.
$(u)_{D^{i}s}=\lfloor u\rfloor_{D^{i}s}$. We prove the result by induction on
$|u|_{D^{i}s}$. If $|u|_{D^{i}s}=1$, then the result holds. Assume
that $|u|_{D^{i}s}>1$, $ \lfloor u\rfloor_{D^{i}s}=b_{(n)}\lfloor u_1\rfloor_{D^{i}s},  \
b\in B, \  i,n\geq 0. $ By induction,
$\lfloor u_1\rfloor_{D^{i}s}=\sum\limits_{j}\beta_{j}[v_j]_{D^{l_j}s_j} $. By
Lemma \ref{l3.3}, we can get the result.
\ \ \ \  $\square$

\subsection{Key lemmas}

The following lemmas play a key role in the proof of the
Composition-Diamond lemma, see Theorem \ref{cdl}.

\begin{lemma}\label{l3.6}
Let $S$ be a set of monic polynomials, closed under the compositions of inclusion
and intersection, and $s_{1}, s_{2}\in S$. If
$w=\overline{[u_{1}]_{D^{i_{1}}s_{1}}}=\overline{[u_{2}]_{D^{i_{2}}s_{2}}}$,
 then
\begin{equation}\label{*}
h:=[u_{1}]_{D^{i_{1}}s_{1}}-[u_{2}]_{D^{i_{2}}s_{2}}
=\sum_t\beta_{_t}[v_{_{t}}]_{_{D^{l_t}s_{_{t}}}}
\end{equation}
where each $[v_{_{t}}]_{_{D^{l_t}s_{_{t}}}}$ is a normal $S$-word  and
$\overline{[v_{_{t}}]_{_{D^{l_t}s_{_{t}}}}}<w$.
\end{lemma}
\noindent{\bf Proof.} If $|u_{1}|_{D^{i_{1}}s_{1}}=1$ or
 $|u_{2}|_{D^{i_{2}}s_{2}}=1$, the result follows from Lemma \ref{l3.3}.
We may assume that $[u_{t}]_{D^{i_{t}}s_{t}}=
[a_{t(n_{t})}D^{i_{t}}s_{t}]$, where $a_{t}\in T $ is $D$-free,
$t=1,2$. There are two cases to consider.

Case 1. Suppose that the subword $\overline{s_{1}}$ of $w$ contains
 $\overline{s_{2}}$ as a subword. Then
 $\overline{s_{1}}=a_{(n_{2})}\overline{s_{2}}D^{l}$,
 $a_{1(n_{1})}a=a_{2}$, $l=i_{2}-i_{1}$, $a\in T$ is
 $D$-free, We have
\begin{align*}
h&=[a_{1(n_{1})}D^{i_{1}}s_{1}]-[a_{2(n_{2})}D^{i_{2}}s_{2}]
=[a_{1(n_{1})}D^{i_{1}}s_{1}]-[a_{1(n_{1})}a_{(n_{2})}D^{i_{2}}s_{2}] \\
&=[a_{1(n_{1})}D^{i_{1}}s_{1}]-[a_{1(n_{1})}D^{i_{1}}[a_{(n_{2})}D^{l}s_{2}]]+\sum_k\alpha_{_k}[q_{_{k}}]_{_{D^{j_k}s_{_{k}}}}\\
&=[a_{1(n_{1})}D^{i_{1}}(s_{1}, s_{2})_{w_1}]+\sum_k\alpha_{_k}[q_{_{k}}]_{_{D^{j_k}s_{_{k}}}},
\end{align*}
where $\overline{[q_{_{k}}]_{_{D^{j_k}s_{_{k}}}}}<\overline{[a_{1(n_{1})}a_{(n_{2})}D^{i_{2}}s_{2}]}=w$,
$(s_{1}, s_{2})_{w_1}$ is the composition of inclusion. So $\overline{D^{i_{1}}(s_{1}, s_{2})_{w_1}}<w_1 D^{i_{1}}$.
By Lemma \ref{l3.3} (ii), we can get the result.

Case 2. Suppose $\overline{s_{1}}$ and $\overline{s_{2}}$ have a nonempty
 intersection as a subword of $w$.
We may assume that  $\overline{s_{1}}D^{l}=a_{(n_{2})}\overline{s_{2}}$,
$a_{1(n_{1})}a=a_{2}$, $l=i_{1}-i_{2}$, $a\in T$ is $D$-free.
Then, similar to Case 1, we can get the result. \ \ \ \  $\square$

\ \

Let $S$ be a set of monic polynomials. Denote
$$
Irr(S)=\{[u]\mid u\in U, u\neq \overline{[v]_{D^{i}s}}\ \
\mbox{for any normal } S\mbox{-word } [v]_{D^{i}s}\}.
$$

\begin{lemma}\label{l3.7}
For any
$0\neq f\in mod_{C(B,N)}(Y)$, $f$ has a presentation
$$
f=\sum_{i}\alpha_{i}[u_{i}]+\sum_{j}\beta_{j}[v_{j}]_{D^{l_{j}}s_{j}},
$$
where each $[u_{i}]\in Irr(S)$, $\overline{[v_{j}]_{D^{l_{j}}s_{j}}}
\leq \bar f$, $\overline{[u_{i}]}\leq \bar f$ and
$s_{j}\in S$.
\end{lemma}
\noindent{\bf Proof.}
If $\bar f=\overline{[u]_{D^{i}s}}$ for some normal $S$-word $[u]_{D^{i}s}$, then let
$f_{1}=f-\alpha[u]_{D^{i}s}$, where $\alpha$ is the leading coefficient of $f$.
If  $[\bar f]\in Irr(S)$, then
let $f_{1}=f-\alpha[\bar f]$.
In both cases, we have $\bar f_1<\bar f$. Thus, the result follows by induction on
$\bar f$. \ \ \ \  $\square$

\subsection{Composition-Diamond lemma for associative conformal modules}

\begin{theorem} (Composition-Diamond lemma for associative conformal modules)\label{cdl}
Let $S\subset mod_{C(B,N)}(Y)$ be a set of monic polynomials, $<$ the ordering on
 $U$ as before and $subm(S)$ the submodule of  $mod_{C(B,N)}(Y)$ generated by $S$.
 Then the following statements are equivalent.
\begin{enumerate}
\item[(i)]\ $S$ is a Gr\"{o}bner-Shirshov basis in $mod_{C(B,N)}(Y)$.
\item[(ii)]\ If $0\neq f\in subm(S)$, then $\bar f
=\overline{[u]_{D^{i}s}}$ for some normal $S$-word $[u]_{D^{i}s}$.
\item[(iii)]\ $Irr(S)=\{[u]\mid u\in U, u\neq \overline{[v]_{D^{i}s}}\ \
\mbox{for any normal } S\mbox{-word } [v]_{D^{i}s}\}$
is a $\mathbf{k}$-basis of the factor module $mod_{C(B,N)}(Y|S):=mod_{C(B,N)}(Y)/subm(S)$.
\end{enumerate}
\end{theorem}
\noindent{\bf Proof.} $(i)\Rightarrow (ii)$. \ Let $S$ be a
Gr\"{o}bner-Shirshov basis and $0\neq f\in subm(S)$. Then, by Lemma
\ref{ll3.4}, we have
$f=\sum\limits_{j}\alpha_{j}[u_{j}]_{D^{l_{j}}s_{j}}$ where each
$\alpha_j\in \mathbf{k}, \ s_j\in S, \ l_{j}\geq 0$. Let
$$
w_j=\overline{[u_{j}]_{D^{l_{j}}s_{j}}}, \ w_1=w_2=\cdots=w_t>w_{t+1}\geq\cdots.
$$
We will use induction on $(w_1,t)$ to prove that
$\bar f=\overline{[u]_{D^{i}s}}$ for some $s\in S$, $i\geq 0$, where $(w_1,t)<(w_1',t')$ lexicographically.

If $w_1=\bar f$, then
the result holds. Assume that $w_1>\bar f, t\geq 2$. We have
$w_1=\overline{[u_{1}]_{D^{l_{1}}s_{1}}}=\overline{[u_{2}]_{D^{l_{2}}s_{2}}}$.
By Lemma \ref{l3.6},
\begin{eqnarray*}
\alpha_{_1}[u_{_{1}}]_{_{D^{l_1}s_{_{1}}}}+\alpha_{_2}[u_{_{2}}]_{_{D^{l_2}s_{_{2}}}}
&=&(\alpha_{_1}+\alpha_{_2})[u_{_{1}}]_{_{D^{l_1}s_{_{1}}}}
+\alpha_{_2}([u_{_{2}}]_{_{D^{l_2}s_{_{2}}}}-[u_{_{1}}]_{_{D^{l_1}s_{_{1}}}})\\
&=&(\alpha_{_1}+\alpha_{_2})[u_{_{1}}]_{_{D^{l_1}s_{_{1}}}}
+\sum\alpha_{_2}\beta_{_t}[v_{_{t}}]_{_{D^{j_t}s_{_{t}}}},
\end{eqnarray*}
where $[v_{_{t}}]_{_{D^{j_t}s_{_{t}}}}$ is a normal $S$-word for each $t$ and
$\overline{[v_{_{t}}]_{_{D^{j_t}s_{_{t}}}}}<w_1$.
Now the result follows from the induction whenever
$\alpha_1+\alpha_2\neq 0$, or $t>2$, or both $\alpha_1+\alpha_2= 0$
and $t=2$. This shows (ii).

$(ii)\Rightarrow (iii)$. By Lemma \ref{l3.7}, $Irr(S)$ generates
$mod_{C(B,N)}(Y|S)$ as a linear space. Suppose that
$\sum\limits_{i}\alpha_i[u_i]=0$ in $mod_{C(B,N)}(Y|S)$, where
$0\neq\alpha_i\in \mathbf{k}$, $[u_i]\in {Irr(S)}$. It means that
$\sum\limits_{i}\alpha_i[u_i]\in{subm(S)}$ in $mod_{C(B,N)}(Y)$.
Then $\overline{\sum\limits_{i}\alpha_i[u_i]}=u_j$ for some $j$.
Since $[u_j]\in{Irr(S)}$, it contradicts (ii).

$(iii)\Rightarrow(i)$.
Suppose that $h$ is any composition of polynomial in $S$. Then $h\in subm(S)$.
Since $Irr(S)$ is a $\mathbf{k}$-basis for $mod_{C(B,N)}(Y)/subm(S)$ and by Lemma \ref{l3.7},
we have $h\equiv 0\ \ mod(S)$,
i.e.
 $S$ is a Gr\"{o}bner-Shirshov basis.
\ \ \ \  $\square$

\ \

Let $C=C(B,N|S)$ be an arbitrary associative conformal algebra. We may assume that $S$ is a Gr\"{o}bner-Shirshov basis  in $C(B,N)$ (see \cite{nll}). Then by Theorem \ref{pro1}, $mod_{C(B,N)}(Y|R)=mod_{C(B,N|S)}(Y)$ is the free $C$-module generated by $Y$, where $R=\{s_{(m)}[u]\mid s\in S, u\in U, m\geq 0\}\subseteq mod_{C(B,N)}(Y)$. Then it is important to find a Gr\"{o}bner-Shirshov basis for $subm(R)$ in $mod_{C(B,N)}(Y)$.

Note that for a $D$-free monic subset $S$ in $C(B,N)$, $S$ is a Gr\"{o}bner-Shirshov basis in $C(B,N)$ if $S$  is closed under the compositions of
inclusion, intersection and left multiplication in $C(B,N)$, see \cite{nll}.

\begin{lemma}\label{l4.1}
Let $S$ be a $D$-free subset of $C(B,N)$,
$R=\{s_{(m)}[u]\mid s\in S, u\in U, m\geq 0\}\subseteq mod_{C(B,N)}(Y)$ and
$R_1=\{s_{(m)}[u]\mid s\in S, u\in U, 0\leq m<N\}$.
Then in $mod_{C(B,N)}(Y)$, $subm(R)=subm(R_{1})$.
\end{lemma}
\noindent{\bf Proof.}  \  We only need to show that
$R\subseteq subm(R_1)$. Since $S$ is $D$-free,
it is clearly that $s_{(p)}b=0, \ s_{(p)}y=0$ for any
$s\in S, \ b\in B, \ y\in Y, \ p\geq N$.
We will prove $s_{(m)}[u]\in span_{\mathbf{k}}R_1$
for any $m\geq N$.

If $[u]=D^{i}y$, we have
$$
s_{(N)}D^{i}y=-\sum_{t\geq 1}(-1)^{t}
\binom{i}{t}\frac{N!}{(N-t)!}s_{(N-t)}D^{i-t}y,
$$
so $s_{(N)}D^{i}y\in span_{\mathbf{k}}R_1$.

Let $|u|>1$ and $[u]=b_{1(n_{1})}[u_{1}]$. Then
$$
s_{(m)}[u]=s_{(m)}(b_{1(n_{1})}[u_{1}])=
-\sum_{t\geq 1}(-1)^{t}\binom{m}{t}s_{(m-t)}(b_{1(n_{1}+t)}[u_{1}]).
$$
By induction on $m$, we complete the proof.\ \ \ \  $\square$

\begin{proposition}\label{p1}
Let the notation be as in Lemma \ref{l4.1} and $S$  a $D$-free Gr\"{o}bner-Shirshov basis in $C(B,N)$.
Then $R_{1}$ is a Gr\"{o}bner-Shirshov basis in $mod_{C(B,N)}(Y)$.
Moreover,
$$
Irr(R_1)=\{[a_{(n)}D^{i}y]\mid [a]\in Irr(S),\ a\in T\ \mbox{with}\ a\ D\mbox{-free},\ 0\leq n<N,\ i\geq 0, y\in Y\},
$$
is a $\mathbf{k}$-basis for the module $mod_{C(B,N)}(Y|R)$, where
$$
Irr(S)=\{[u]\mid u\in T,\ u \mbox{ is } D\mbox{-free and}\ u\neq \overline{[v]_{D^{i}s}}\ \
\mbox{for any normal } s\mbox{-word } [v]_{D^{i}s}\}.
$$
\end{proposition}
\noindent{\bf Proof.}  \ By Lemma \ref{l4.1}, we have $subm(R)=subm(R_1)$.

For any $f=s_{1(m_{1})}[u_{1}]$, $g=s_{2(m_{2})}[u_{2}]\in R_{1}$,
we assume that $[u_{j}]=[a_{j(n_j)}D^{i_j}y_j], \ j=1,2$. Hence,
\begin{align*}
&f-[s_{1(m_{1})}a_1]_{(n_1)}D^{i_1}y_1\\
&=f-\sum_{t\geq 0}(-1)^{t}\binom{m_1}{t}s_{1(m_{1}-t)}([a_1]_{(n_1+t)}D^{i_1}y_1) \\
&=f-s_{1(m_{1})}([a_1]_{(n_1)}D^{i_1}y_1)-\sum_{t\geq 1}(-1)^{t}\binom{m_1}{t}s_{1(m_{1}-t)}([a_1]_{(n_1+t)}D^{i_1}y_1)\\
&=\sum_{k}\alpha_k h_k,
\end{align*}
where $h_k\in R_1$ and $\overline{h_k}<\bar{f}$.

Similarly, we have
$g-[s_{2(m_{2})}a_2]_{(n_2)}D^{i_2}y_2=\sum_{t}\gamma_t q_t$,
where $q_t\in R_1$ and $\overline{q_t}<\bar{f}$.

If $f,g$ have composition, then $n_1=n_2=:n$,
$y_1=y_2=:y$. There are three cases to consider.

Case 1. $w=\bar f=a_{(p)}\bar gD^{k}$ where $a\in T$ is $D$-free and $k\geq 0$. Then
$\bar s_{1(m_{1})}a_1=a_{(p)}\bar s_{2(m_{2})}a_2=:c$,
and $i_1=i_2+k$. Since $S$ is a $D$-free Gr\"{o}bner-Shirshov basis in $C(B,N)$, we have
$$
[s_{1(m_{1})}a_1]-[a_{(p)}s_{2(m_{2})}a_2]=\sum_{j'}\beta_{_j'}[v_{_{j'}}]_{_{D^{l_{j'}}s_{_{j'}}}},
$$
where each $[v_{_{j'}}]_{_{D^{l_{j'}}s_{_{j'}}}}$ is a normal $S$-word and
$\overline{[v_{_{j'}}]_{_{D^{l_{j'}}s_{_{j'}}}}}<c$.
Hence,
\begin{align*}
f-[a_{(p)}D^{l}g]
&= [s_{1(m_{1})}a_1]_{(n)}D^{i_1}y-[a_{(p)}D^{l}([s_{2(m_{2})}a_2]_{(n)}D^{i_2}y)]\\
&\ \ \  +\sum_{k}\alpha_k h_k-\sum_{t}\gamma_t[a_{(p)}D^{l}q_t] \\
&= [s_{1(m_{1})}a_1]_{(n)}D^{i_1}y-[a_{(p)}([s_{2(m_{2})}a_2]_{(n)}D^{l+i_2}y)]\\
&\ \ \ -\alpha [a_{(p)}([s_{2(m_{2})}a_2]_{(n-l)}D^{i_2}y)]
+\sum_{k}\alpha_k h_k-\sum_{t}\gamma_t[a_{(p)}D^{l}q_t]  \\
&= ([s_{1(m_{1})}a_1]-[a_{(p)}s_{2(m_{2})}a_2])_{(n)}D^{i_1}y-\alpha [a_{(p)}([s_{2(m_{2})}a_2]_{(n-l)}D^{i_2}y)]\\
&\ \ \ +\sum_{k}\alpha_k h_k-\sum_{t}\gamma_t[a_{(p)}D^{l}q_t]\\
&\equiv 0\  mod (R_1),
\end{align*}
where $\alpha=(-1)^{l}\frac{n!}{(n-l)!}$.

Case 2.  $w=\bar f D^{i}=a_{(n)}\bar g$ where $a\in T$ is $D$-free and $i\geq 0$. Similar to Case 1, we have $D^{i}f-[a_{(n)}g]\equiv 0\  mod (R_1)$.

Case 3. For any $b\in B,\ n\geq N$, we will prove
$b_{(n)}f\equiv 0 \  mod (R_1)$. If $n=N$,
then by Lemma \ref{l4.1},
\begin{align*}
b_{(N)}f&=b_{(N)}(s_{1(m_{1})}[u_{1}])=(b_{(N)}s_1)_{(m_1)}[u_1]-
\sum_{t\geq 1}(-1)^{t}\binom{N}{t}b_{(N-t)}(s_{1(m_{1}+t)}[u_{1}]).
\end{align*}
Since $S$ is a $D$-free Gr\"{o}bner-Shirshov basis, we have
$b_{(N)}s_1=\sum\limits_j \beta_j [a_{j(n_j)}s_{j(q_j)}c_j]$,
where $a_j,\ c_j\in T$ are $D$-free. So
\begin{align*}
b_{(N)}f
&\equiv \sum\limits_j \beta_j [a_{j(n_j)}s_{j(q_j)}c_j]_{(m_1)}[u_1]
\equiv \sum\limits_j \beta_j [a_{j(n_j)}s_{j(q_j)}c_{j(m_1)}u_1] \equiv 0  \ mod(R_1).
\end{align*}

\noindent Assume that $n>N$. By induction, we have
$
b_{(n)}f\equiv (b_{(n)}s_1)_{(m_1)}[u_1]\equiv 0\ mod(R_1).
$

Therefore, $R_1$ is a Gr\"{o}bner-Shirshov basis in $mod_{C(B,N)}(Y)$. Now, by Theorem \ref{cdl}, the
set $R_1$ is a $\mathbf{k}$-basis for the module $mod_{C(B,N)}(Y|R)$.
\ \ \ \  $\square$

\ \

\noindent{\bf Remark:} \
The condition that  $S$ is $D$-free in Proposition \ref{p1}
is essential. For example, let $C=C(a,N=2|S)$ be an associative conformal algebra, where
$$
S=\{a_{(1)}a-a_{(0)}Da,\  [a_{(0)}a_{(1)}a],\ [a_{(1)}a_{(0)}a],\ [a_{(0)}a_{(0)}a],\ [a_{(1)}a_{(1)}a]\}.
$$
It is easy to check that $S$ is a Gr\"{o}bner-Shirshov basis in
$C(a,N=2)$. Let $Y$ be a well-ordered set,  $mod_{C(B,N=2)}(Y)$ the
double free conformal module and $R=\{s_{(m)}[u]\mid s\in S, u\in
U, m\geq 0\}$, where $U=\{a_{(n_{1})}\cdots a_{(n_{k})}D^{i}y\mid
y\in Y, 0\leq n_{j}<2, 1\leq j\leq k, i, k\geq 0\}$. Since
$(a_{(1)}a-a_{(0)}Da)_{(2)}y=a_{(0)}(Da_{(2)}y)=-2a_{(0)}(a_{(1)}y)$,
$a_{(0)}(a_{(1)}y)\in subm(R)$. Thus
$(a_{(1)}a-a_{(0)}Da)_{(2)}y=0\in mod_{C(B,N=2)}(Y|R)$. Noting that
$a_{(0)}a\in Irr(S)$,  the set $\{[a_{(n)}D^{i}y]\mid [a]\in
Irr(S) \}$ isn't a $\mathbf{k}$-basis of $mod_{C(B,N=2)}(Y|R)$.

\section{Applications}
\subsection{Conformal modules over universal enveloping conformal algebra}

Let $L$ be a Lie conformal algebra which is a free
$\mathbf{k}[D]$-module with a well-ordered $\mathbf{k}[D]$-basis $B=\{a_i\mid
i\in I\}$ and a uniform bounded locality $N(a_i,a_j)\leq N$ for all
$i, j\in I$. Let the multiplication table of $L$ on $B$ be
$$
a_i{_{[n]}}a_j=\Sigma_{_{t\in I}}\alpha_{nijt}a_t,\ \ \alpha_{nijt}\in
\mathbf{k}[D],\ i\geq j,\ i,j\in I,\ n<N.
$$
Then by $\mathcal{U}(L)$, the universal enveloping associative conformal algebra of $L$ with respective to $B$ and $N$,
one means the following associative conformal algebra, see \cite{BFK04},
$$
\mathcal{U}(L)=C(B,N\mid a_i{_{(n)}}a_j-\{a_j{_{(n)}}a_i\}-\lfloor a_i{_{[n]}}a_j\rfloor,\ i\geq j,\ i,j\in I, n<N)
$$
where $\{a_j{_{(n)}}a_i\}=\sum_{k\geq
0}(-1)^{n+k}\frac{1}{k!}D^{k}(a_j{_{(n+k)}}a_i)$ and $\lfloor
a_i{_{[n]}}a_j\rfloor=\Sigma_{_{t\in I}}\alpha_{nijt}a_t$.

Let $\mathbb{C}$ be the complex field,
$$
Vir=C_{Lie}(v,\ N=2\mid v_{[0]}v-Dv, v_{[1]}v-2v)
$$
be the Lie conformal algebra over  $\mathbb{C}$
($Vir$ is called the Virasoro conformal algebra), see \cite{BFK04}, and
$$
\mathcal{U}(Vir)=C(v,\ N=2\mid v_{(1)}v-v) 
$$
the universal enveloping associative conformal algebra of $Vir$.

\begin{example}\label{ex1}
Let $\Delta\in\{0, 1\}$ and $\alpha\in\mathbb{C}$. Let
$$
M(\Delta, \alpha)= mod_{C(v; N=2)}(y| R\cup Q),
$$
where
$$
R=\{(v_{(1)}v-v)_{(m)}[v_{(n_1)}v\cdots v _{(n_k)}D^{i}y]\mid n_j\in \{0,1\},\ 1\leq j\leq k,\ m, k, i\geq 0 \}
$$
and $Q=\{f_1, f_2\}$ where $f_1=v_{(0)}y-(D+\alpha)y,\ f_2=v_{(1)}y-\Delta y$.

Then $Q=\{f_1, f_2\}$ is a Gr\"{o}bner-Shirshov basis for
$ M(\Delta, \alpha)$.
 So, by Theorem \ref{cdl}, the set $Irr(Q)=\{D^{i}y\mid i\geq 0\}$
is a $\mathbb{C}$-basis of $ M(\Delta, \alpha)$.
It follows that
$M(\Delta, \alpha)=\mathbb{C}[D]y$ is a ${\mathcal{U}(Vir)}$-module, called the Virasoro conformal module, see \cite{SJVK, Kac96, ARe}.

\end{example}
\noindent{\bf Proof.}
Let $R_1=\{(v_{(1)}v-v)_{(m)}[v_{(n_1)}v\cdots v _{(n_k)}D^{i}y]\mid m,n_j\in \{0,1\},\ 1\leq j\leq k,\ k, i\geq 0 \}$.
By Lemma 4.2 in  \cite{BFK04}, $\{v_{(1)}v-v\}$ is a $D$-free Gr\"{o}bner-Shirshov basis in $C(v; N=2)$. Then by Lemma \ref{l4.1}, $subm(R)=subm(R_1)$ in $mod_{C(v; N)}(y)$.
So $M(\Delta, \alpha)= mod_{C(v; N)}(y| R\cup Q)=mod_{C(v; N)}(y| R_1\cup Q)$.

If $i>0$, then
\begin{align*}
v_{(m)}D^{i}y&=\sum\limits_{t\geq 0}\binom{m}{t}\frac{i!}{(i-t)!}D^{i-t}(v_{(m-t)}y) \\
&= \left\{
\begin{aligned}
& D^{i}(v_{(0)}y),\ \ m=0,   \\
&D^{i}(v_{(1)}y)+iD^{i-1}(v_{(0)}y),\ \ m=1
\end{aligned}
\right. \\
&\equiv\left\{
\begin{aligned}
& (D+\alpha)D^{i}y,\ \ m=0,   \\
&\Delta D^{i}y+i(D+\alpha)D^{i-1}y,\ \ m=1
\end{aligned}
\right. \\
&\equiv \left\{
\begin{aligned}
&D^{i+1}y+\alpha D^{i}y\ \ mod(Q),\ \ m=0,   \\
&(\Delta+i)D^{i}y+i\alpha D^{i-1}y\ \ mod(Q),\ \ m=1.
\end{aligned}
\right.
\end{align*}
So $[v_{(n_1)}v\cdots v _{(n_k)}D^{i}y]\equiv \sum\limits_{l\geq 0}\beta_{l}D^{l}y\ mod(Q)$, where
$n_j\in \{0,1\},\ 1\leq j\leq k,\ k, i\geq 0$ and $\beta_{l}\in \mathbb{C}$.
Denote $s=v_{(1)}v-v$. Let $h=s_{(m)}[v_{(n_1)}v\cdots v _{(n_k)}D^{i}y]\in R_1$. Then we can get $h\equiv 0\ mod(Q)$ following
$s_{(n)}D^{k}y\equiv 0\ mod(Q)$ for any $k\geq 0, n\in\{0,1\}$.

Now,
\begin{eqnarray*}
v_{(2)}Dy&=&\sum\limits_{t\geq 0}\binom{2}{t}\frac{1}{(1-t)!}D^{1-t}(v_{(2-t)}y)\equiv 2\Delta y\ mod(Q),\\
v_{(2)}D^{i}y&=&\sum\limits_{t\geq 0}\binom{2}{t}\frac{i!}{(i-t)!}D^{i-t}(v_{(2-t)}y)=2iD^{i-1}(v_{(1)}y)+i(i-1)D^{i-2}(v_{(0)}y)\\
&\equiv&(i^{2}-i+2\Delta i)D^{i-1}y+\alpha(i^{2}-i)D^{i-2}y\ mod(Q),\ i\geq 2,\\
s_{(0)}y&=&(v_{(1)}v)_{(0)}y-v_{(0)}y
=v_{(1)}(v_{(0)}y)-v_{(0)}(v_{(1)}y)-v_{(0)}y \\
&\equiv& v_{(1)}(D+\alpha)y-v_{(0)}\Delta y-(D+\alpha)y \\
&\equiv& v_{(1)}Dy-(\Delta+1)Dy-\alpha y \\
&\equiv& 0 \ mod(Q),
\end{eqnarray*}
\begin{eqnarray*}
s_{(0)}D^{i}y&=&(v_{(1)}v)_{(0)}D^{i}y-v_{(0)}D^{i}y
=v_{(1)}(v_{(0)}D^{i}y)-v_{(0)}(v_{(1)}D^{i}y)-v_{(0)}D^{i}y \\
&\equiv& v_{(1)}(D^{i+1}y+\alpha D^{i}y)-v_{(0)}((\Delta+i)D^{i}y+i\alpha D^{i-1}y)-D^{i+1}y-\alpha D^{i}y \\
&\equiv& (\Delta+i+1)D^{i+1}y+(i+1)\alpha D^{i}y+\alpha(\Delta+i)D^{i}y+i\alpha^{2}D^{i-1}y \\
&&\ \ \ \  -(\Delta+i)(D^{i+1}y+\alpha D^{i}y)-i\alpha D^{i}y-i\alpha^{2}D^{i-1}y-D^{i+1}y-\alpha D^{i}y \\
&\equiv& 0 \ mod(Q), \ \ \ \ \  i>0,\\
s_{(1)}D^{i}y&=&\sum\limits_{t\geq 0}\binom{1}{t}\frac{i!}{(i-t)!}D^{i-t}(s_{(1-t)}y)
=\sum\limits_{t\geq 0}\binom{1}{t}\frac{i!}{(i-t)!}D^{i-t}((v_{(1)}v)_{(1-t)}y-v_{(1-t)}y)\\
&=&\sum\limits_{t\geq 0}\binom{1}{t}\frac{i!}{(i-t)!}D^{i-t}(v_{(1)}(v_{(1-t)}y)-v_{(0)}(v_{(2-t)}y)-v_{(1-t)}y)\\
&\equiv& \left\{
\begin{aligned}
&v_{(1)}(v_{(1)}y)-v_{(1)}y ,\ \ i=0,   \\
&iD^{i-1}(v_{(1)}(v_{(0)}y)-v_{(0)}(v_{(1)}y)-v_{(0)}y)+D^{i}(v_{(1)}(v_{(1)}y)-v_{(1)}y),\ \ i>0
\end{aligned}
\right.\\
&\equiv& \left\{
\begin{aligned}
&(\Delta^{2}-\Delta)y ,\ \ i=0,   \\
&iD^{i-1}(v_{(1)}(D+\alpha)y-v_{(0)}\Delta y-(D+\alpha)y)+(\Delta^{2}-\Delta)D^{i}y,\ \ i>0
\end{aligned}
\right. \\
&\equiv& \left\{
\begin{aligned}
&(\Delta^{2}-\Delta)y ,\ \ i=0,   \\
&iD^{i-1}(v_{(1)}Dy-(\Delta+1)Dy-\alpha y)+(\Delta^{2}-\Delta)D^{i}y,\ \ i>0
\end{aligned}
\right. \\
&\equiv& (\Delta^{2}-\Delta)D^{i}y \\
&\equiv& 0 \ mod(Q).
\end{eqnarray*}

For any $n\geq 2$, we have
\begin{eqnarray*}
v_{(n)}f_1&=&v_{(n)}(v_{(0)}y)-v_{(n)}(D+\alpha)y \\
&=&-\sum\limits_{t\geq 1}(-1)^{t}\binom{n}{t}v_{(n-t)}(v_{(t)}y)-v_{(n)}Dy-\alpha v_{(n)}y \\
&\equiv& nv_{(n-1)}(v_{(1)}y)-nv_{(n-1)}y\\
&\equiv& \left\{
\begin{aligned}
& 0,\ \ \ \ n>2,   \\
&2(\Delta^{2}-\Delta)y,\ \ \ n=2,
\end{aligned}
\right.\\
&\equiv& 0 \ mod(Q),\\
v_{(n)}f_2&=&v_{(n)}(v_{(1)}y)-v_{(n)}\Delta y
=-\sum\limits_{t\geq 1}(-1)^{t}\binom{n}{t}v_{(n-t)}(v_{(1+t)}y)-\Delta v_{(n)}y \\
&\equiv& 0 \ mod(Q).
\end{eqnarray*}
From this it follows that $subm(Q)=subm(R_1\cup Q)$ and  all left multiplication compositions in $Q$
are trivial modulo $Q$.

Then $Q$ is closed under the left
multiplication composition. Since $Q$ has no composition of
inclusion and intersection, $Q$ is a Gr\"{o}bner-Shirshov basis of
$M(\Delta, \alpha)$.

Now, by Theorem  \ref{cdl} and Proposition \ref{p3}, the results follow. \ \ \ \  $\square$

\begin{example}
 Module over the semidirect product of Virasoro conformal algebra and current algebra.

Let $(\mathfrak{g},[\ ])$ be a Lie algebra over $\mathbb{C}$ with a well-ordered $\mathbb{C}$-basis $\{a_i\}_{i\in I}$ and
$Cur(\mathfrak{g})$ be the current algebra over $\mathfrak{g}$, where
$$
Cur(\mathfrak{g})=C_{Lie}(\{a_i\}_{i\in I},\ N=1\mid
a_{i[0]}a_j=[a_ia_j], \ i,j\in I).
$$
The semidirect product of $ Vir$ and $Cur(\mathfrak{g})$ is
$$
Vir\oplus Cur(\mathfrak{g})=C_{Lie}(\{v\}\cup\{a_i\}_{i\in I},\ N=2\mid S ),
$$
where
$$
S=\{v_{[0]}v-Dv, v_{[1]}v-2v, v_{[0]}a_i-Da_i, v_{[1]}a_i-a_i, a_{i[0]}a_j-[a_ia_j], a_{i[1]}a_j,\ i>j,\ i,j\in I \}.
$$
Then, see section 4.4 in \cite{BFK04},
$$
\mathcal{U}(Vir\oplus Cur(\mathfrak{g}))=C(\{v\}\cup\{a_i\}_{i\in I},\ N=2\mid S^{(-)}),
$$
where  $S^{(-)}$ consists of
\begin{align*}
&s_1=a_{i(0)}a_j-a_{j(0)}a_i-[a_ia_j],\ i>j,\\
&s_2=a_{i(1)}a_j,\ i>j,\\
&s_3=v_{(1)}v-v,  \\
&s_4=v_{(0)}a_i+a_{i(1)}Dv-2a_{i(0)}v-Da_i,\\
&s_5=v_{(1)}a_i+a_{i(1)}v-a_i, \\
&s_6=v_{(0)}(a_{j(0)}a_i)-a_{j(0)}(v_{(0)}a_i),\ i>j,\\
&s_7=v_{(0)}(a_{i(0)}v)-a_{i(0)}(v_{(0)}v), \\
&s_8=v_{(0)}(a_{i(1)}v)-a_{i(1)}(v_{(0)}v)+a_{i(0)}v.
\end{align*}

Let $V$ be a $\mathfrak{g}$-module with a $\mathbb{C}$-basis $Y$ and
$\Delta\in\{0, 1\},\ \alpha\in \mathbb{C}$. Let
$$
M(\Delta, \alpha,V)=
mod_{C(\{v\}\cup\{a_i\}_{i\in I}, N)}(Y| R\cup Q),
$$
where
\begin{align*}
R=\{s_{(m)}[c_1{_{(n_1)}}\cdots c_k{_{(n_k)}}D^{t}y]\mid\  &s\in S^{(-)},\ n_j\in \{0,1\},\ c_j\in\{v\}\cup\{a_i\}_{i\in I},\\
&y\in Y,\ 1\leq j\leq k,\ k, m,t\geq 0 \}
\end{align*}
and $Q=\{f_{1y}, f_{2y}, f_{3iy}, f_{4iy}\mid  i\in I, y\in Y\}$ where
$$
f_{1y}=v_{(0)}y-(D+\alpha)y,\ f_{2y}= v_{(1)}y-\Delta y,\ f_{3iy}=a_{i}{_{(0)}}y-a_iy,\ f_{4iy}=a_{i}{_{(1)}}y,\  i\in I, y\in Y.
$$

Then $Q$ is a Gr\"{o}bner-Shirshov basis for $ M(\Delta, \alpha,V)$. Moreover,
$M(\Delta, \alpha,
V)={\mathbb{C}[D]}\!Y$ is a  $Vir\oplus Cur(\mathfrak{g})$-module, see \cite{SJVK, Kac96}.
\end{example}
\noindent{\bf Proof.} Let $m\geq 0, k\geq 1$. Then we have, $mod(Q)$,
\begin{eqnarray*}
v_{(m)}D^{k}y&=&\sum\limits_{t\geq 0}\binom{m}{t}\frac{k!}{(k-t)!}D^{k-t}(v_{(m-t)}y) \\
&\equiv& \left\{
\begin{aligned}
&\frac{k!}{(k-m)!}D^{k-m}(v_{(0)}y)+m\frac{k!}{(k-m+1)!}D^{k-m+1}(v_{(1)}y),\ \ k\geq m,   \\
&\binom{m}{k}v_{(m-k)}y,\ \ k<m
\end{aligned}
\right.\\
&\equiv& \left\{
\begin{aligned}
&\frac{k!}{(k-m)!}(D+\alpha)D^{k-m}y+m\Delta\frac{k!}{(k-m+1)!}D^{k-m+1}y,\ \ m\leq k,   \\
&(k+1)\Delta y,\ \ m=k+1,\\
&0,\ \ \ \ m>k+1,
\end{aligned}
\right.\\
a_{i(m)}D^{k}y&=&\sum\limits_{t\geq 0}\binom{m}{t}\frac{k!}{(k-t)!}D^{k-t}(a_{i(m-t)}y)\\
&\equiv& \left\{
\begin{aligned}
& 0,\ \ \ \ \ 1\leq k<m,   \\
&\frac{k!}{(k-m)!}D^{k-m}(a_{i}y),\ \ \ \  k\geq m\geq 0.
\end{aligned}
\right.
\end{eqnarray*}
So $[c_1{_{(n_1)}}\cdots c_k{_{(n_k)}}D^{i}y]\equiv \sum\limits_{l\geq 0}\beta_{l}D^{l}y\ mod(Q)$, where
$n_j\in \{0,1\},\ c_j\in\{v\}\cup\{a_i\}_{i\in I},\ 1\leq j\leq k,\ i,k\geq 0$ and $\beta_{l}\in \mathbb{C}$.
Let $h=s_{(m)}[c_1{_{(n_1)}}\cdots c_k{_{(n_k)}}D^{i}y]\in R$. Then we can get $h\equiv 0\ mod(Q)$ following
$s_{(n)}D^{k}y\equiv 0\ mod(Q)$ for any $s\in S^{(-)}, n,k\geq 0$.

Let $m\geq 0$. For $s_1, s_2, s_3, s_6$, we have, $mod(Q)$,
\begin{eqnarray*}
s_{1(m)}D^{k}y&=&(a_{i(0)}a_j-a_{j(0)}a_i-[a_ia_j])_{(m)}D^{k}y \\
&=&a_{i(0)}(a_{j(m)}D^{k}y)-a_{j(0)}(a_{i(m)}D^{k}y)-[a_ia_j]_{(m)}D^{k}y \\
&\equiv& \left\{
\begin{aligned}
& 0,\ \  1\leq k<m,   \\
&\frac{k!}{(k-m)!}D^{k-m}(a_i(a_jy))-(a_j(a_iy))-[a_{i}a_j]y),\ \   k\geq m\geq 0
\end{aligned}
\right. \\
&\equiv& 0,\\
s_{2(m)}D^{k}y&=&\sum\limits_{t\geq 0}\binom{m}{t}\frac{k!}{(k-t)!}D^{k-t}(s_{2(m-t)}y)\\
&=&\sum\limits_{t\geq 0}\binom{m}{t}\frac{k!}{(k-t)!}D^{k-t}((a_{i(1)}a_j)_{(m-t)}y) \\
&=&\sum\limits_{t\geq 0}\binom{m}{t}\frac{k!}{(k-t)!}D^{k-t}((a_{i(1)}(a_{j(m-t)}y)-(a_{i(0)}(a_{j(m-t+1)}y))\\
&\equiv& 0,\\
s_{3(0)}D^{k}y&=&D^{k}((v_{(1)}v-v)_{(0)}y)=D^{k}(v_{(1)}(v_{(0)}y)-v_{(0)}(v_{(1)}y)-v_{(0)}y) \\
&\equiv& D^{k}(v_{(1)}(D+\alpha)y-v_{(0)}\Delta y-(D+\alpha)y) \\
&\equiv& 0,
\end{eqnarray*}
\begin{eqnarray*}
s_{3(m)}D^{k}y&=&\sum\limits_{t\geq 0}\binom{m}{t}\frac{k!}{(k-t)!}D^{k-t}(s_{3(m-t)}y)\\
&=&\sum\limits_{t\geq 0}\binom{m}{t}\frac{k!}{(k-t)!}D^{k-t}((v_{(1)}v)_{(m-t)}y-v_{(m-t)}y) \\
&=&\sum\limits_{t\geq 0}\binom{m}{t}\frac{k!}{(k-t)!}D^{k-t}(v_{(1)}(v_{(m-t)}y)-v_{(0)}(v_{(m-t+1)}y)-v_{(m-t)}y) \\
&\equiv&\left\{
\begin{aligned}
& 0,\ \  m>k+1\geq 1,   \\
&m\frac{k!}{(k+1-m)!}D^{k+1-m}(\Delta^{2}-\Delta)y,\ \   1\leq m\leq k+1,
\end{aligned}
\right.\\
&\equiv& 0,\\
s_{6(m)}D^{k}y&=&\sum\limits_{t\geq 0}\binom{m}{t}\frac{k!}{(k-t)!}D^{k-t}((v_{(0)}(a_{j(0)}a_i)-a_{j(0)}(v_{(0)}a_i))_{(m-t)}y) \\
&=&\sum\limits_{t\geq 0}\binom{m}{t}\frac{k!}{(k-t)!}D^{k-t}(v_{(0)}(a_{j(0)}(a_{i(m-t)}y))-a_{j(0)}(v_{(0)}(a_{i(m-t)}y))) \\
&\equiv&\left\{
\begin{aligned}
& 0,\ \  m>k\geq 0,   \\
&\frac{k!}{(k-m)!}D^{k-m}(v_{(0)}(a_{j(0)}\{a_iy\})-a_{j(0)}(v_{(0)}\{a_iy\})),\ \   0\leq m\leq k
\end{aligned}
\right. \\
&\equiv& 0.
\end{eqnarray*}

For $s_4, s_5, s_7, s_8$ and $m\geq 0$, we have
\begin{eqnarray*}
s_{4(m)}D^{k}y&=&\sum\limits_{t\geq 0}\binom{m}{t}\frac{k!}{(k-t)!}D^{k-t}((v_{(0)}a_i+a_{i(1)}Dv-2a_{i(0)}v-Da_i)_{(m-t)}y) \\
&=&\sum\limits_{t\geq 0}\binom{m}{t}\frac{k!}{(k-t)!}D^{k-t}(v_{(0)}(a_{i(m-t)}y)+a_{i(1)}(Dv_{(m-t)}y)-a_{i(0)}(Dv_{(m-t+1)}y) \\
&& \ \ \ \ \ \ \ \ \ \ \ \ \ \ \ \ \ \ \ \ \ \ \ \ \ \ \ \ \ \ \ \ \ \ -2a_{i(0)}(v_{(m-t)}y)-Da_{i(m-t)}y) \\
&=& \sum\limits_{t\geq 0}\binom{m}{t}\frac{k!}{(k-t)!}D^{k-t}(v_{(0)}(a_{i(m-t)}y)+a_{i(1)}(Dv_{(m-t)}y) \\
&& \ \ \ \ \ \ \ \ \ \ \ \ \ \ \ \ \ \ \ \ \ \ \ \ \ \ \ \ \ \ \ \ \ \
+(m-t+1)a_{i(0)}(v_{(m-t)}y)-Da_{i(m-t)}y),\\
s_{5(m)}D^{k}y&=&\sum\limits_{t\geq 0}\binom{m}{t}\frac{k!}{(k-t)!}D^{k-t}((v_{(1)}a_i+a_{i(1)}v-a_i)_{(m-t)}y) \\
&=&\sum\limits_{t\geq 0}\binom{m}{t}\frac{k!}{(k-t)!}D^{k-t}(v_{(1)}(a_{i(m-t)}y)-v_{(0)}(a_{i(m-t+1)}y)+a_{i(1)}(v_{(m-t)}y) \\
&&\ \ \ \ \ \ \ \ \ \ \ \ \ \ \ \ \ \ \ \ \ \ \ \ \ \ \ \ \ \ \ \ \ \ -a_{i(0)}(v_{(m-t+1)}y)-a_{i(m-t)}y),\\
s_{7(m)}D^{k}y&=&\sum\limits_{t\geq 0}\binom{m}{t}\frac{k!}{(k-t)!}D^{k-t}((v_{(0)}(a_{i(0)}v)-a_{i(0)}(v_{(0)}v))_{(m-t)}y) \\
&=&\sum\limits_{t\geq 0}\binom{m}{t}\frac{k!}{(k-t)!}D^{k-t}(v_{(0)}(a_{i(0)}(v_{(m-t)}y))-a_{i(0)}(v_{(0)}(v_{(m-t)}y))),
\end{eqnarray*}
\begin{eqnarray*}
s_{8(m)}D^{k}y&=&\sum\limits_{t\geq 0}\binom{m}{t}\frac{k!}{(k-t)!}D^{k-t}((v_{(0)}(a_{i(1)}v)-a_{i(1)}(v_{(0)}v)+a_{i(0)}v)_{(m-t)}y) \\
&=&\sum\limits_{t\geq 0}\binom{m}{t}\frac{k!}{(k-t)!}D^{k-t}(v_{(0)}(a_{i(1)}(v_{(m-t)}y))-v_{(0)}(a_{i(0)}(v_{(m-t+1)}y)) \\
&&\ \ \ \ \ \ \ \ \ \ \ \ \ \ \ \ -a_{i(1)}(v_{(0)}(v_{(m-t)}y))+a_{i(0)}(v_{(0)}(v_{(m-t+1)}y))+a_{i(0)}(v_{(m-t)}y)).
\end{eqnarray*}
There are two cases to consider.

Case 1.  $m>k\geq 0$. Then, $mod(Q)$,
\begin{eqnarray*}
s_{4(m)}D^{k}y&\equiv & \binom{m}{k}k!(a_{i(1)}(Dv_{(m-k)}y)+(m-k-1)a_{i(0)}(v_{(m-k)}y)-Da_{i(m-k)}y) \\
&\equiv& \left\{
\begin{aligned}
& 0,\ \  m-k>1,   \\
&m!(-a_{i(1)}(v_{(0)}y)+a_{i(0)}y ),\ \   m-k=1
\end{aligned}
\right. \\
&\equiv& 0,\\
s_{5(m)}D^{k}y&\equiv& \binom{m}{k}(v_{(1)}(a_{i(m-k)}y)-v_{(0)}(a_{i(m-k+1)}y)+a_{i(1)}(v_{(m-k)}y) \\
&\equiv& 0,\\
s_{7(m)}D^{k}y&=&\binom{m}{k}(v_{(0)}(a_{i(0)}(v_{(m-k)}y))-a_{i(0)}(v_{(0)}(v_{(m-k)}y)))\\
&\equiv& \left\{
\begin{aligned}
& 0,\ \  m-k>1,   \\
&m(v_{(0)}(a_{i(0)}\Delta y)-a_{i(0)}(v_{(0)}\Delta y)),\ \   m-k=1
\end{aligned}
\right.\\
&\equiv& \left\{
\begin{aligned}
& 0,\ \  m-k>1,   \\
&m(v_{(0)}\Delta (a_iy)-a_{i(0)}\Delta(D+\alpha)y),\ \   m-k=1
\end{aligned}
\right. \\
&\equiv& 0,\\
s_{8(m)}D^{k}y&\equiv& \binom{m}{k}(v_{(0)}(a_{i(1)}(v_{(m-k)}y))-a_{i(1)}(v_{(0)}(v_{(m-k)}y))+a_{i(0)}(v_{(m-k)}y))\\
&\equiv& \left\{
\begin{aligned}
& 0,\ \  m-k>1,   \\
&m(v_{(0)}(a_{i(1)}\Delta y)-a_{i(1)}(v_{(0)}\Delta y)+a_{i(0)}\Delta y),\ \   m-k=1
\end{aligned}
\right. \\
&\equiv&  \left\{
\begin{aligned}
& 0,\ \  m-k>1,   \\
&m(v_{(0)}(-a_{i(1)}\Delta(D+\alpha)y+\Delta(a_iy),\ \   m-k=1
\end{aligned}
\right.\\
&\equiv& 0.
\end{eqnarray*}

Case 2. $0\leq m\leq k$. Then, $mod(Q)$,
\begin{eqnarray*}
s_{4(m)}D^{k}y&=&\frac{k!}{(k-m)!}D^{k-m}(v_{(0)}(a_{i(0)}y)-a_{i(0)}(v_{(0)}y)) \\
&&+m\frac{k!}{(k+1-m)!}D^{k+1-m}(a_{i(1)}(Dv_{(1)}y)+a_{i(0)}y)\\
&\equiv& \frac{k!}{(k-m)!}D^{k-m}(v_{(0)}(a_iy)-a_{i(0)}(D+\alpha)y) \\
&&\ \ \  +m\frac{k!}{(k+1-m)!}D^{k+1-m}(-a_{i(1)}(D+\alpha)y+a_{i(0)}y) \\
&\equiv& 0,
\end{eqnarray*}
\begin{eqnarray*}
s_{5(m)}D^{k}y&\equiv& \frac{k!}{(k-m)!}D^{k-m}(v_{(1)}(a_{i(0)}y)+a_{i(1)}(v_{(0)}y)-a_{i(0)}(v_{(1)}y)-a_{i(0)}y )\\
&\equiv& 0,\\
s_{7(m)}D^{k}y&=&\frac{k!}{(k-m)!}D^{k-m}(v_{(0)}(a_{i(0)}(v_{(0)}y))-a_{i(0)}(v_{(0)}(v_{(0)}y)))\\
&&\ \ \ \ +m\frac{k!}{(k-m+1)!}D^{k-m+1}(v_{(0)}(a_{i(0)}(v_{(1)}y))-a_{i(0)}(v_{(0)}(v_{(1)}y)))\\
&\equiv& \frac{k!}{(k-m)!}D^{k-m}(v_{(0)}(a_{i(0)}(D+\alpha)y)-a_{i(0)}(v_{(0)}(D+\alpha)y))\\
&&\ \ \ \  +m\frac{k!}{(k-m+1)!}D^{k-m+1}(v_{(0)}(a_{i(0)}\Delta y)-a_{i(0)}(v_{(0)}\Delta y))\\
&\equiv& \frac{k!}{(k-m)!}D^{k-m}(v_{(0)}(D+\alpha)(a_iy))-a_{i(0)}(D+\alpha)^{2}y)\\
&&\ \ \ \  +m\frac{k!}{(k-m+1)!}D^{k-m+1}(v_{(0)}\Delta (a_iy)-a_{i(0)}\Delta(D+\alpha)y)\\
&\equiv& 0,\\
s_{8(m)}D^{k}y&\equiv& \frac{k!}{(k-m)!}D^{k-m}(v_{(0)}(a_{i(1)}(v_{(0)}y))-v_{(0)}(a_{i(0)}(v_{(1)}y))-a_{i(1)}(v_{(0)}(v_{(0)}y)) \\
&&\ \ \ \ \ \ \ \ \ \ \ \ \ \ \ \  \ \ \ \ \ \ \ \ +a_{i(0)}(v_{(0)}(v_{(1)}y))+a_{i(0)}(v_{(0)}y))\\
&&\ \  +m\frac{k!}{(k-m+1)!}D^{k-m+1}(v_{(0)}(a_{i(1)}(v_{(1)}y))-a_{i(1)}(v_{(0)}(v_{(1)}y))+a_{i(0)}(v_{(1)}y))\\
&\equiv& 0.
\end{eqnarray*}

Let $n\geq 2$ and $i,j\in I$. Then we have, $mod(Q)$,
\begin{eqnarray*}
v_{(n)}f_{1y}&\equiv &\left\{
\begin{aligned}
& 0,\ \ \ n>2,   \\
&2(\Delta^{2}-\Delta)y,\ \ \ n=2,
\end{aligned}
\right.
\equiv 0,\\
v_{(n)}f_{2y}&=&v_{(n)}(v_{(1)}y)-v_{(n)}\Delta y\equiv 0,
\\
v_{(n)}f_{3iy}&=&v_{(n)}(a_{i(0)}y-a_iy)\equiv \sum\limits_{k\geq 1}(-1)^{k+1}\binom{n}{k}v_{(n-k)}(a_{i(k)}y)
\equiv 0,\\
v_{(n)}f_{4iy}&=&v_{(n)}(a_{i(1)}y)=\sum\limits_{k\geq 1}(-1)^{k+1}\binom{n}{k}v_{(n-k)}(a_{i(k+1)}y)
\equiv 0 ,
\end{eqnarray*}
\begin{eqnarray*}
a_{j(n)}f_{1y}&=&a_{j(n)}(v_{(0)}y-(D+\alpha)y)=a_{j(n)}(v_{(0)}y)-a_{j(n)}(D+\alpha)y) \\
&\equiv& \sum\limits_{k\geq 1}(-1)^{k+1}\binom{n}{k}a_{j(n-k)}(v_{(k)}y)-na_{j(n-1)}y \\
&\equiv& na_{j(n-1)}(v_{(1)}y)
\equiv 0,\\
a_{j(n)}f_{2y}&=&a_{j(n)}(v_{(1)}y-\Delta y)
\equiv \sum\limits_{k\geq 1}(-1)^{k+1}\binom{n}{k}a_{j(n-k)}(v_{(1+k)}y)\equiv 0,\\
a_{j(n)}f_{3iy}&=&a_{j(n)}(a_{i(0)}y-a_iy)
\equiv\sum\limits_{k\geq 1}(-1)^{k+1}\binom{n}{k}a_{j(n-k)}(a_{i(k)}y)\equiv 0,\\
a_{j(n)}f_{4iy}&=&a_{j(n)}(a_{i(1)}y)\equiv 0.
\end{eqnarray*}

This shows that $subm(Q)=subm(R\cup Q)$ and  all left multiplication compositions in $Q$
are trivial modulo $Q$.
 Since $Q$ has no composition of inclusion and
intersection, $Q$ is a Gr\"{o}bner-Shirshov basis of $M(\Delta,
\alpha,V)$. 

Now, the results follow from Theorem  \ref{cdl} and Proposition \ref{p3}. \ \ \ \  $\square$

\end{document}